\documentclass[12pt,letterpaper,3p]{elsarticle}

\usepackage{layout}
\usepackage[lined,boxed,commentsnumbered]{algorithm2e}
\usepackage{makeidx}
\usepackage{amsthm}
\usepackage{amsmath}
\usepackage{amscd}
\usepackage{mathrsfs}
\usepackage{bbold}
\usepackage{geometry} 
\usepackage{enumitem}
\usepackage{setspace}
\usepackage{titlesec}
\usepackage{titletoc}
\usepackage[parfill]{parskip}    
\usepackage{graphicx}
\usepackage{amssymb}
\usepackage{subfigure}
\usepackage{tikz}
\usetikzlibrary{arrows}
\usepackage{tikz-cd}
\usepackage{algorithm2e}
\usepackage[all]{xy}

\usepackage{epstopdf}
\DeclareGraphicsExtensions{.eps}

\swapnumbers

\newcommand{\blob}{
 \rule[.0ex]{1ex}{1ex}
}

\swapnumbers

\theoremstyle{definition}
\newtheorem{theorem}{Theorem}[section]

\newcommand{\dive}{\operatorname{div}}

\newcommand{\diag}{\operatorname{diag}}

\newcommand{\supp}{\operatorname{supp}}

\newcommand{\inte}{\operatorname{int}}
\newcommand{\ind}{\operatorname{ind}}

\makeatletter
\renewcommand*\env@matrix[1][*\c@MaxMatrixCols c]{%
  \hskip -\arraycolsep
  \let\@ifnextchar\new@ifnextchar
  \array{#1}}
\makeatother

\makeatletter
\def\ps@pprintTitle{%
  \let\@oddhead\@empty
  \let\@evenhead\@empty
  \def\@oddfoot{\reset@font\hfil\thepage\hfil}
  \let\@evenfoot\@oddfoot
}
\makeatother

\begin{document}

\begin{frontmatter}




\title{A Preconditioner Based on Low-Rank Approximation of Schur Complements}


\author[pg]{P.~Gatto}
\ead{paolo.gatto@epfl.ch}

\author[pg]{J.S.~Hesthaven}
\ead{jan.hesthaven@epfl.ch}

\address[pg]{Mathematics Institute of Computational Science and Engineering \\
 \'{E}cole Polytechnique F\'{e}d\'{e}rale de Lausanne (EPFL) \\
 MA C2 652 (Batiment MA), Station 8 \\ 
 CH-1015 Lausanne, Switzerland}

\begin{abstract}
We introduce a preconditioner based on a hierarchical low-rank compression scheme of Schur complements. The construction is inspired by standard nested dissection, and relies on the assumption that the Schur complements can be approximated, to high precision, by Hierarchically-Semi-Separable matrices. We build the preconditioner as an approximate $LDM^t$ factorization of a given matrix $A$, and no knowledge of $A$ in assembled form is required by the construction. The $LDM^t$ factorization is amenable to fast inversion, and the action of the inverse can be determined fast as well. We investigate the behavior of the preconditioner in the context of DG finite element approximations of elliptic and hyperbolic problems.
\end{abstract}

\begin{keyword}
Preconditioned GMRES \sep Interpolative Decomposition \sep Indefinite operators

\end{keyword}

\end{frontmatter}

\section{Introduction}\label{sec:0}

This work rests on the observation that, for a large class of problems, the dense Schur complement matrices that arise in the nested dissection method are rank-structured, see, e.g., \cite{borm2010approximation, xia2009superfast}. In the case of well-behaved elliptic problems, this property is a consequence of the rapid decay of the underlying Green function. To the contrary, in the case of wave-propagation problems, the same argument does not apply, and the reasons behind the rank-structure of the Schur complements remain poorly understood. Nevertheless, by exploiting this property, approximate matrix decompositions can be constructed cheaply, and turn out to be excellent preconditioners.

Linear systems that arise from finite element discretizations of wave propagation phenomena are typically poorly conditioned and highly indefinite. Consequently, it is both vital and challenging to construct effective preconditioners. Standard multigrid methods generally fail to carry the oscillations at the wavelength scale onto the coarse grids. Incomplete LU-decompositions require to assemble the global matrix, are fairly expensive to compute, and still lead to a number of iterations that is frequency-dependent. The lack of effective preconditioning techniques has lead to the use of (sparse) direct solvers. Fast methods, such as the fast-multipole-method, are confined to problems for which the governing PDE's can be reformulated as boundary integral equations (BIE's). The linear systems resulting from the discretization of the BIE's are dense, as opposed to sparse, and often well-conditioned. Over the last twenty years, a number of methods has been developed for their efficient solution. They are based on a rigorous understanding of the physics of the problem, along with sophisticated analytical arguments that rely on asymptotic expansions of special functions. As a result, they have limited applicability, e.g., variable coefficients problems are out of reach, and their efficient implementation remains quite challenging.

As of today, because of the lack of robust preconditioners, discretizations of wave propagation problems have eluded the reach of fast iterative solvers. In this work, we introduce a preconditioner whose construction is completely general, is parallel in nature, and fits modern computational architectures. Results obtained in the context of Discontinuous Garlerkin (DG) finite element approximations show that the behavior of the preconditioner, measured as the number of GMRES iterations, is independent of both the mesh size and the order of approximation.

The paper is organized as follows. In Section \ref{sec:1} we review the concept of Hierarchically-Semi-Separable matrices, and describe how a matrix can be reduced to such form within linear complexity. In Section \ref{sec:2} we develop analytical arguments to justify the low-rank nature of the Schur complements. The construction of the preconditioner is described in full details in Section \ref{sec:3}, and numerical examples that substantiate its robustness are reported in Section \ref{sec:4}. Finally, in Section \ref{sec:5}, we draw conclusions from this work, and point towards future directions of research.

\section{Rank-Structured Matrices}\label{sec:1}

In general terms, a matrix is rank-structured if the ranks of its off-diagonal blocks, obtained through a recursive partitioning of the index vector, are small when compared to the matrix size. A number of rank-structured formats have been described in the literature, e.g., the $\mathcal{H}$- and $\mathcal{H}^2$-matrices introduced in \cite{hackbusch1999sparse,hackbusch2000sparse}. We focus on Hierarchically-Semi-Separable (HSS) matrices. Their off-diagonal blocks admit low-rank decompositions into factors satisfying certain recursion relations that make such matrices inexpensive to store and manipulate. In fact, if $k$ is the off-diagonal rank of a square HSS matrix of size $N$, then such matrix can be applied to a vector and inverted in 
$\mathrm{O}(Nk)$ and $\mathrm{O}(Nk^2)$ operations, respectively.

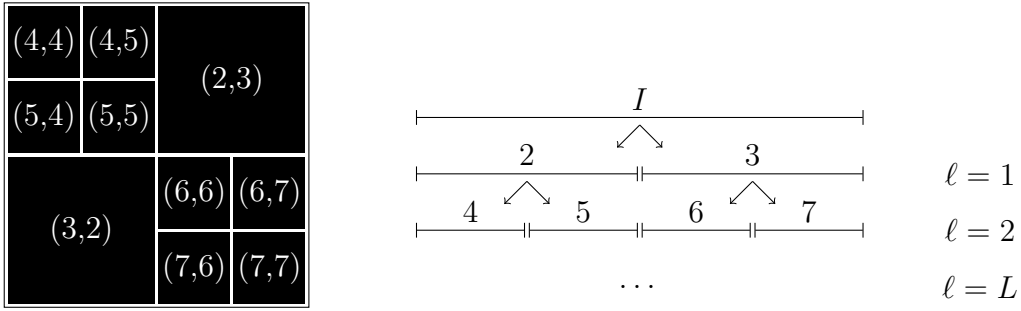
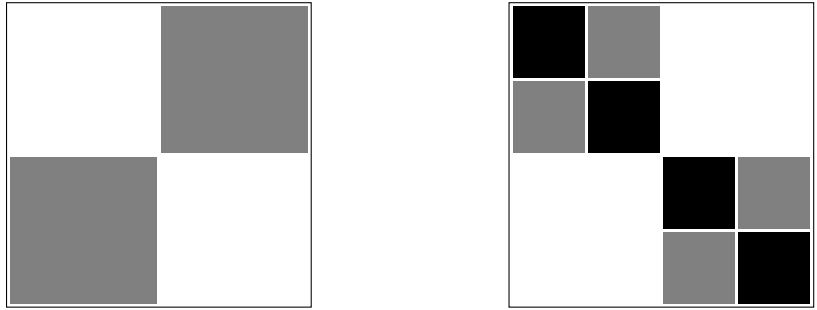
\begin{figure}[t]

\begin{center}
\subfigure[Block indexing of matrix $A$ resulting from partitioning of index vector $I$.\label{fig:0.1}]{
\begin{tikzpicture}[scale = 0.05]

\draw[black] (0,0) rectangle (81,81);

\path[fill=black] (1,1) rectangle (40,40);
\node[white] at(20.5,20.5) {(3,2)};

\path[fill=black] (41,1) rectangle (60,20);
\node[white] at(50.5,10.5) {(7,6)};

\path[fill=black] (61,1) rectangle (80,20);
\node[white] at(70.5,10.5) {(7,7)};

\path[fill=black] (41,21) rectangle (60,40);
\node[white] at(50.5,30.5) {(6,6)};

\path[fill=black] (61,21) rectangle (80,40);
\node[white] at(70.5,30.5) {(6,7)};

\path[fill=black] (41,41) rectangle (80,80);
\node[white] at(60.5,60.5) {(2,3)};

\path[fill=black] (1,41) rectangle (20,60);
\node[white] at(10.5,50.5) {(5,4)};

\path[fill=black] (21,41) rectangle (40,60);
\node[white] at(30.5,50.5) {(5,5)};

\path[fill=black] (1,61) rectangle (20,80);
\node[white] at(10.5,70.5) {(4,4)};

\path[fill=black] (21,61) rectangle (40,80);
\node[white] at(30.5,70.5) {(4,5)};

\end{tikzpicture}
\qquad\quad
\begin{tikzpicture}[scale = 0.05]

\draw[black, |-|] (0,45) -- (119,45);
\node at (59.5,50) {$I$};
\draw[black,->] (59.5,43)--(53.5,37);
\draw[black,->] (59.5,43)--(65.5,37);

\draw[black, |-|] (0,30) -- (59,30);
\draw[black, |-|] (60,30) -- (119,30);

\draw[black,->] (29.5,28)--(23.5,22);
\draw[black,->] (29.5,28)--(35.5,22);

\draw[black,->] (89.5,28)--(83.5,22);
\draw[black,->] (89.5,28)--(95.5,22);

\node at (150,30) {$\ell = 1$};
\node at (29.5,35) {$2$};
\node at (89.5,35) {$3$};

\draw[black, |-|] (0,15) -- (29,15);
\draw[black, |-|] (30,15) -- (59,15);
\draw[black, |-|] (60,15) -- (89,15);
\draw[black, |-|] (90,15) -- (119,15);
\node at (150,15) {$\ell = 2$};
\node at (14.5,20) {$4$};
\node at (44.5,20) {$5$};
\node at (74.5,20) {$6$};
\node at (104.5,20) {$7$};

\node at (59.5,0) {$\dots$};
\node at (150,0) {$\ell = L$};
\end{tikzpicture}
} \\

\vspace{0.3cm}

\subfigure[Low-rank blocks, indicated in gray, for level $\ell= 1$ (left), as opposed to low-rank and full-rank blocks, indicated in black, for leaf-level $\ell=2$ (right).\label{fig:0.2}]{
\begin{tikzpicture}[scale = 0.05]

\draw[black] (0,0) rectangle (81,81);

\path[fill=gray] (1,1) rectangle (40,40);
\path[fill=gray] (41,41) rectangle (80,80);

\end{tikzpicture}
 \qquad \qquad \qquad
\begin{tikzpicture}[scale = 0.05]

\draw[black] (0,0) rectangle (81,81);

\path[fill=gray] (41,1) rectangle (60,20);

\path[fill=black] (61,1) rectangle (80,20);

\path[fill=black] (41,21) rectangle (60,40);

\path[fill=gray] (61,21) rectangle (80,40);

\path[fill=gray] (1,41) rectangle (20,60);

\path[fill=black] (21,41) rectangle (40,60);

\path[fill=black] (1,61) rectangle (20,80);

\path[fill=gray] (21,61) rectangle (40,80);

\end{tikzpicture}
}
\end{center}

\caption{Rank-structure of matrix $A$.  \label{fig:0}}
\end{figure}

To make the discussion precise, let $A$ be a square matrix of size $N= 2^L$ and partition its index vector $I = (1, \ldots, 2^Lm)$ recursively through a binary tree with $L$ levels, see Figure \ref{fig:0.1}. Following a well-establish terminology, we call leaf nodes those nodes that belong to the finest level, namely $\ell = L$. Nodes $\sigma$, $\tau$ on level $\ell$ form a sibling pair if they share a common ancestor on level $\ell -1$. Matrix $A$ is called an ``$\mathcal{S}$-matrix'' or a ``semi-separable matrix'' if there exists an integer $k$ such that, for every sibling pair $\{ \sigma ,\tau\}$ of the tree, the off-diagonal blocks $A(\sigma,\tau)$\footnote{Throughout the paper, the \textsc{Matlab}\textsuperscript{\textregistered}-like notation $A(\sigma,\tau)$ indicates the restriction of $A$ to row-index vector $\sigma$ and column-index vector $\tau$. } and $A(\tau,\sigma)$ have (numerical) ranks equal to $k$, see Figure \ref{fig:0.2}.

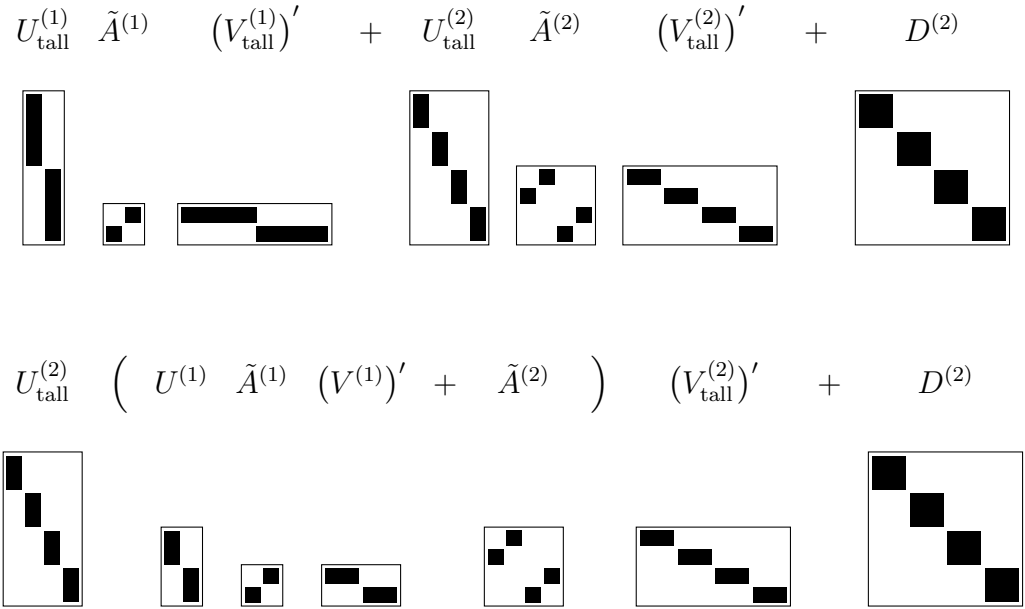
\begin{figure}[t]

\[
\begin{array}{cccc cccc c}
U^{(1)}_{\text{tall}} & \tilde{A}^{(1)} & {\big( V^{(1)}_{\text{tall}} \big)}' & + &
U^{(2)}_{\text{tall}} & \tilde{A}^{(2)} & {\big( V^{(2)}_{\text{tall}} \big)}' & + &
D^{(2)} \\
\\

\begin{tikzpicture}[scale = 0.05]
\draw[black] (0,0) rectangle (11,41);
\path[fill=black] (1,21) rectangle (5,40);
\path[fill=black] (6,1) rectangle (10,20);
\end{tikzpicture}
& 

\begin{tikzpicture}[scale = 0.05]
\draw[black] (0,0) rectangle (11,11);
\path[fill=black] (1,1) rectangle (5,5);
\path[fill=black] (6,6) rectangle (10,10);
\end{tikzpicture}
&

\begin{tikzpicture}[scale = 0.05]
\draw[black] (0,0) rectangle (41,11);
\path[fill=black] (1,6) rectangle (21,10);
\path[fill=black] (21,1) rectangle (40,5);
\end{tikzpicture}
&

&
\begin{tikzpicture}[scale = 0.05]
\draw[black] (0,0) rectangle (21,41);
\path[fill=black] (1,31) rectangle (5,40);
\path[fill=black] (6,21) rectangle (10,30);
\path[fill=black] (11,11) rectangle (15,20);
\path[fill=black] (16,1) rectangle (20,10);
\end{tikzpicture} 
&

\begin{tikzpicture}[scale = 0.05]
\draw[black] (0,0) rectangle (21,21);
\path[fill=black] (11,1) rectangle (15,5);
\path[fill=black] (16,6) rectangle (20,10);

\path[fill=black] (1,11) rectangle (5,15);
\path[fill=black] (6,16) rectangle (10,20);
\end{tikzpicture}
&

\begin{tikzpicture}[scale = 0.05]

\draw[black] (0,0) rectangle (41,21);
\path[fill=black] (31,1) rectangle (40,5);
\path[fill=black] (21,6) rectangle (30,10);
\path[fill=black] (11,11) rectangle (20,15);
\path[fill=black] (1,16) rectangle (10,20);
\end{tikzpicture}

& 
&
\begin{tikzpicture}[scale = 0.05]
\draw[black] (0,0) rectangle (41,41);
\path[fill=black] (31,1) rectangle (40,10);
\path[fill=black] (21,11) rectangle (30,20);
\path[fill=black] (11,21) rectangle (20,30);
\path[fill=black] (1,31) rectangle (10,40);
\end{tikzpicture}

\end{array}
\] \\

\begin{equation*}
\begin{array}{cccc cccc ccc}
U_{\text{tall}}^{(2)} & \Big( & U^{(1)} & \tilde{A}^{(1)} &   \big(V^{(1)}\big)'  & + &  \tilde{A}^{(2)} & \Big) & \big(V_{\text{tall}}^{(2)}\big)'  & + & D^{(2)} \\
\\

\begin{tikzpicture}[scale = 0.05]
\draw[black] (0,0) rectangle (21,41);
\path[fill=black] (1,31) rectangle (5,40);
\path[fill=black] (6,21) rectangle (10,30);
\path[fill=black] (11,11) rectangle (15,20);
\path[fill=black] (16,1) rectangle (20,10);
\end{tikzpicture} & &

\begin{tikzpicture}[scale = 0.05]
\draw[black] (0,0) rectangle (11,21);
\path[fill=black] (1,11) rectangle (5,20);
\path[fill=black] (6,1) rectangle (10,10);
\end{tikzpicture} &

\begin{tikzpicture}[scale = 0.05]
\draw[black] (0,0) rectangle (11,11);
\path[fill=black] (1,1) rectangle (5,5);
\path[fill=black] (6,6) rectangle (10,10);
\end{tikzpicture} &

\begin{tikzpicture}[scale = 0.05]
\draw[black] (0,0) rectangle (21,11);
\path[fill=black] (11,1) rectangle (20,5);
\path[fill=black] (1,6) rectangle (10,10);
\end{tikzpicture} & &

\begin{tikzpicture}[scale = 0.05]
\draw[black] (0,0) rectangle (21,21);
\path[fill=black] (11,1) rectangle (15,5);
\path[fill=black] (16,6) rectangle (20,10);

\path[fill=black] (1,11) rectangle (5,15);
\path[fill=black] (6,16) rectangle (10,20);
\end{tikzpicture} & &

\begin{tikzpicture}[scale = 0.05]

\draw[black] (0,0) rectangle (41,21);
\path[fill=black] (31,1) rectangle (40,5);
\path[fill=black] (21,6) rectangle (30,10);
\path[fill=black] (11,11) rectangle (20,15);
\path[fill=black] (1,16) rectangle (10,20);
\end{tikzpicture} & &

\begin{tikzpicture}[scale = 0.05]
\draw[black] (0,0) rectangle (41,41);
\path[fill=black] (31,1) rectangle (40,10);
\path[fill=black] (21,11) rectangle (30,20);
\path[fill=black] (11,21) rectangle (20,30);
\path[fill=black] (1,31) rectangle (10,40);
\end{tikzpicture}

\end{array}
\end{equation*}

\caption{Factorization of a semi-separable matrix (top), as opposed to the factorization of a hierarchically-semi-separable matrix (bottom), for $L = 2$. For an HSS matrix, at the non-leaf levels, the ``tall'' basis matrices $U_{\text{tall}}^{(\cdot)}$ and $V_{\text{tall}}^{(\cdot)}$ are replaced by their ``short'' counterparts  $U^{(\cdot)}$ and $V^{(\cdot)}$.}\label{fig:1}

\end{figure}

For each sibling pair $\{\sigma,\tau\}$ on level $\ell$, the corresponding off-diagonal blocks of a semi-separable matrix can be factored as:
\begin{equation*}
A(\sigma,\tau)   =  U_{\sigma}^{\text{tall}} \,  \tilde{A}_{\sigma,\tau}  \, {\big( V_{\tau}^{\text{tall}} \big)}'  \qquad ,  \qquad
A(\tau,\sigma)  =  U_{\tau}^{\text{tall}} \, \tilde{A}_{\tau,\sigma} \, {\big( V_{\sigma}^{\text{tall}} \big)}' 
\end{equation*}
where square matrices $\tilde{A}_{\sigma,\tau}$ and $\tilde{A}_{\tau,\sigma}$ have size $k$, which is independent of $\ell$. Let us define the block-diagonal matrices:
\begin{align*}
D^{(L)} & = \diag\{ D_\sigma = A(\sigma,\sigma) \:\: : \:\:  \tau \text{ belongs to level } L \} \\
U^{(\ell)}_{\text{tall}} & = \diag\{ U_\sigma^\text{tall} \:\: : \:\:  \sigma \text{ belongs to level } \ell \} && \text{for } \ell = 1, \ldots, L \\
V^{(\ell)}_{\text{tall}} & = \diag\{ V_\sigma^\text{tall} \:\: : \:\:  \sigma \text{ belongs to level } \ell \} && \text{for } \ell = 1, \ldots, L
\end{align*}
and let $\tilde{A}^{(\ell)}$ be comprised of all blocks $\tilde{A}_{\sigma,\tau}$ such that $\{ \sigma,\tau\}$ is a sibling pair on level $\ell$. A semi-separable matrix can be factorized as follows:
\begin{equation}\label{eq:10}
A = U^{(1)}_{\text{tall}} \, \tilde{A}^{(1)} \, {\big( V^{(1)}_{\text{tall}} \big)}' + \cdots + U^{(L)}_{\text{tall}} \, \tilde{A}^{(L)} \, {\big( V^{(L)}_{\text{tall}} \big)}' + D^{(L)}
\end{equation}
A pictorial description of the factorization is shown in Figure \ref{fig:1}.

Semi-separable matrices require to store and manipulate ``tall'' $U$- and $V$-factors at all tree-levels. This undesirable property can be circumvented by imposing one additional condition. A semi-separable matrix is said to be Hierarchically-Semi-Separable (HSS) if it satisfies:
\[
U^{\text{tall}}_\sigma =
\begin{pmatrix}
U^{\text{tall}}_\mu &  \\
& U_\nu^{\text{tall}}
\end{pmatrix} U_\sigma
\qquad ; \qquad
V^{\text{tall}}_\sigma =
\begin{pmatrix}
V_\mu^{\text{tall}} &  \\
& V_\nu^{\text{tall}}
\end{pmatrix} V_\sigma
\]
for every $\sigma$ belonging to tree-level $\ell = 1,\ldots, L-1$, with descending sibling pair $\{ \mu, \nu\}$. This is essentially a condition about nesting of column and row spaces. By collecting matrices $U_\sigma$ and $V_\sigma$ into block-diagonal matrices
\begin{align*}
U^{(\ell)} & = \diag\{ U_\sigma \:\: : \:\:  \sigma \text{ belongs to level } \ell \} && \text{for } \ell = 1, \ldots, L-1 \\
V^{(\ell)} & = \diag\{ V_\sigma \:\: : \:\:  \sigma \text{ belongs to level } \ell \} && \text{for } \ell = 1, \ldots, L-1
\end{align*}
the following recursions are obtained:
\[
U^{(\ell)}_{\text{tall}}  = U_\text{tall}^{(L)} \, U^{(L-1)} \cdots U^{(\ell)} \quad , \quad V^{(\ell)}_{\text{tall}}  = V_{\text{tall}}^{(L)} \, V^{(L-1)} \cdots V^{(\ell)}  \qquad \text{for } \ell = 1, \ldots , L-1
\]
Thus, in the case of an HSS matrix, factorization (\ref{eq:10}) simplifies to:
\begin{subequations}
\begin{align*}
A & = U_{\text{tall}}^{(L)} \, \cdots \, U^{(1)} \, \tilde{A}^{(1)} \, {V^{(1)}}' \, \cdots \, {V_{\text{tall}}^{(L)}}' + \cdots + U_{\text{tall}}^{(L)} \, \tilde{A}^{(L)} \, {V_{\text{tall}}^{(L)}}' + D^{(L)} \\
& = U_{\text{tall}}^{(L)} \bigg( U^{(L-1)} \,  \cdots \, U^{(1)} \, \tilde{A}^{(1)} \, {V^{(1)}}' \, \cdots \, {V^{(L-1)}}' + \cdots + \tilde{A}^{(L)} \bigg) {V_{\text{tall}}^{(L)}}' + D^{(L)} \\
& = U_{\text{tall}}^{(L)} \bigg( U^{(L-1)} \Big( U^{(L-2)}  \,  \cdots \, U^{(1)} \, \tilde{A}^{(1)} \, {V^{(1)}}' \, \cdots \, {V^{(L-2)}}' + \cdots + \tilde{A}^{(L-1)} \Big)  V^{(L-1)} \\
& \qquad\qquad\qquad\qquad\qquad\qquad\qquad\qquad\qquad\qquad\qquad\qquad +  \tilde{A}^{(L)} \bigg) {V_{\text{tall}}^{(L)}}' + D^{(L)} \\
& = \cdots
\end{align*}
\end{subequations}
We remark that the assumption that the (numerical) rank of sibling interactions is constant 
across the entire tree is, in practice, replaced by an assumption of boundedness, see Section \ref{sec:3.2} and Section \ref{sec:4} for further discussion.

Rather than HSS matrices \emph{per se}, the construction of  HSS approximants is of practical interest. More precisely, given a matrix $A$ and a tolerance $\varepsilon$, we seek an HSS matrix $A^{\text{(HSS)}}$ such that $\| A -  A^{\text{(HSS)}}\| \le \varepsilon$, for some suitable norm $\| \cdot \|$. If $A$ is an arbitrary square matrix of size $N$, a straightforward approach for computing $A^{\text{(HSS)}}$ yields a $\mathrm{O}(k N^2)$ cost, where $k$ is the HSS-rank. In practice, this is prohibitively expensive. Nevertheless, under moderate assumptions on $A$, it is possible to construct $A^{\text{(HSS)}}$ at the acceptable cost of $\mathrm{O}(Nk^2)$, see \cite{Martinsson:2011:FRA:2340911.2340918}. Thus, provided that $k$ is independent\footnote{As we shall see in Section \ref{sec:4}, this is never the case in practice, and some dependency of $k$ upon $N$ is to be expected.} of $N$, $A$ can be reduced to HSS-form within linear complexity. Let us recall the exact result.

\begin{theorem}\label{th:1}
Let $A$ be an $N \times N$ hierarchically-semi-separable matrix that has HSS-rank $k$. Suppose that:
\begin{enumerate}
\item matrix-vector products $ x \mapsto Ax$ and $ x \mapsto A^tx$ can be evaluated at a cost $T_\text{mult}$;\label{th:item:1}
\item individual entries of $A$ can be evaluated at a cost $T_\text{entry}$.
\end{enumerate}
An HSS factorization of $A$ can be computed in a time proportional to
\[
T_\text{mult} \times 2(k + p) + T_\text{rand} \times N ( k + p ) + T_\text{entry} \times 2 N k + T_\text{flop} \times c N k^2
\]
where $T_\text{rand}$ is the time required to generate a random number, $T_\text{flop}$ is the time required to perform a floating point operation, $p$ is a small oversampling parameter (typically $p=10$), and $c$ is a small constant independent of $N$ or $k$. $\blob$
\end{theorem}

The construction described in \cite{Martinsson:2011:FRA:2340911.2340918} relies heavily on the nesting of the basis of the off-diagonal blocks, and on the compression of such blocks through Interpolative Decompositions (ID), see, e.g., Section 4 and 5 of \cite{chan1992some} for a detailed discussion on interpolative decompositions. In general terms, if $A$ is an $m \times n$ matrix, an ID is a factorization of the form:
 \[
A= A^\text{(skel)}P
 \]
where $ A^\text{(skel)}$ is a ``skeletonization'' of $A$, namely it is an $m \times k$ matrix constructed by selecting $k$ columns of $A$, and $P$ is a well-behaved\footnote{In the present context, by ``well-behaved'', we refer to the fact that no entry of $P$ has absolute value greater than $2$. Let us recall that in general it is possible to obtain an ID where the entries of $P$ are bounded in absolute valued by $1$, although this is an NP-hard problem.} $k \times n$ matrix that, obviously, contains an identity of size $k$. The cost of computing a rank-$k$ ID of $A$ is $\mathrm{O}\big(mkn\log(n)\big)$. As shown in \cite{Woolfe2008}, this cost can be lowered to $\mathrm{O}\big(mn \log (l) + lkn \log(n)\big)$, where $l$ is an integer greater than but close to $k$ (in applications, $l = k + 10$ is typical.) The first term is the cost of obtaining $l$ samples of the range of $A$ \emph{via} an accelerated fast Fourier transform. Consequently, if we additionally assume that $A$ can be applied to a vector within linear complexity, we can drop the first term, and the cost of computing the ID reduces to
\begin{equation}\label{eq:43}
\mathrm{O}\big( lkn \log(n)\big)
\end{equation}
This is the cornerstone to the establish the complexity of the algorithm on which rests the proof of Theorem \ref{th:1}.

\section{Analytical Apparatus}\label{sec:2}

In this section we provide an insight into the rank-structure of the Schur complements that arise in the construction of the preconditioner. To make the discussion precise, let us consider the boundary value problem:
\begin{subequations}\label{pb:0}
\begin{alignat}{3}
\mathcal{L} \, u & = f  \qquad && \text{in } &&\Omega \\
\mathcal{B} \, u & = g && \text{on } &&\Gamma := \partial \Omega
\end{alignat}
\end{subequations}
where $\mathcal{L}$ is a linear second order differential operator, and $\mathcal{B}$ is a linear boundary condition. Whenever the problem is well-posed, there exists a Green's function $G$ such that:
\[
u(x) =\int_\Omega G(x,y) \, f(y) \, dy + \int_\Gamma H(x,y) \, g(y) \, dS(y)
\]
where $H$ is related to the Green's function through the boundary condition $\mathcal{B}$. The right-hand-side of the previous equality defines the so-called solution operator for problem (\ref{pb:0}). A rigorous derivation of the previous integral equation rests upon a generalized Green's formula of the second type, while a characterization of $G$ as the solution of a boundary value problem involves the notion of formal adjoint operator, see, e.g, Chapter 6.7 of \cite{oden2010applied}.

When $\mathcal{L}$ is a uniformly elliptic operator with smooth coefficients and $\Omega$ has a smooth boundary, then the long-range interactions of $G$ are rank-deficient. Recently, the case of non smooth, $L^\infty$-coefficients, and $\Omega$ a bounded Lipschitz domain has been addressed in \cite{bebendorf2003existence}. The authors show that, even in this general setting, the Green's function $G$ can be approximated to high precision by an $\mathcal{H}$-matrix. In the case of homogeneous boundary conditions, this property immediately extends to the solution operator. When non-homogeneous boundary conditions are considered, we postulate that the nature of the solution operator is preserved.

Let $A$ be the stiffness matrix arising from a finite element discretization of problem (\ref{pb:0}). Let us suppose that, up to a permutation we shall omit, we can partition  $A$ as follows, and define the (aggregated) submatrices $A^{(k)}$:
\[
A = 
\begin{pmatrix}
{A^{(1)}}_{ii}	&				& {A^{(1)}}_{ib}	& 				\\
			& {A^{(2)}}_{ii}	&				& {A^{(2)}}_{ib}	\\
{A^{(1)}}_{bi}	&				& {A^{(1)}}_{bb}	& {A^{(1,2)}} 		\\
			& {A^{(2)}}_{bi}	& {A^{(2,1)}}		& {A^{(2)}}_{bb}
\end{pmatrix} \qquad , \qquad
A^{(k)} = \begin{pmatrix}
{A^{(k)}}_{ii} & {A^{(k)}}_{ib} \\
{A^{(k)}}_{bi} &{A^{(k)}}_{bb}
\end{pmatrix} \quad k =1,2
\]
In the partitioning of $A$, we assume the blocks on the diagonal to be square matrices. We proceed to characterize matrix $A^{(k)}$ as the discretization of a boundary value problem.

In the finite element method, each degree of freedom $j$ of the stiffness matrix corresponds to a unique finite element basis function $\varphi_j$. We define the following subdomains of $\Omega$ \emph{via} the supports of the basis functions:
\[
\overline{\Omega_i^{(k)}}  = \cup \, \{ \supp \varphi_j \: : \: j \in \ind ( {A^{(k)}}_{ii}) \} \quad , \quad
\overline{\Omega_b^{(k)}}  = \cup \, \{ \supp \varphi_j \: : \: j \in \ind ( {A^{(k)}}_{bb})  \} 
\]
The notation $\ind (B)$ refers to the (row or column) indices of $A$ that are indices of its submatrix $B$ as well. We also define $\overline{\Omega^{(k)}} = \overline{\Omega_i^{(k)}} \cup \overline{\Omega_b^{(k)}}$. When the previous discretization scheme is applied to the boundary value problem
\begin{subequations}\label{pb:1}
\begin{alignat}{3}
\mathcal{L} \, u & = f  \qquad && \text{in } &&\Omega^{(k)} \\
\mathcal{B} \, u & = g && \text{on } &&\Gamma \cap \partial \Omega^{(k)} \\
u & = 0 && \text{on } &&\inte (\partial \Omega^{(k)} \setminus \Gamma )
\end{alignat}
\end{subequations}
it gives rise to $A^{(k)}$ as a stiffness matrix.

An $LDM^t$ block-factorization of $A^{(k)}$ is realized as follows:
\[
A^{(k)} =
\underset{L}{
\begin{pmatrix}
I 				&  \\
{A^{(k)}}_{ib} \, {{A^{(k)}}_{ii}}^{-1}	& I
\end{pmatrix}
}
\begin{pmatrix}
{A^{(k)}}_{ii}	& \\
 		& S^{(k)}
\end{pmatrix}
\underset{M^t}{
\begin{pmatrix}
I 	& {{A^{(k)}}_{ii}}^{-1} \, {A^{(k)}}_{ib} \\
 	& I
\end{pmatrix}
}
\]
Matrices $L$, $M$ are Gauss transforms, and the Schur complement $S^{(\cdot)}$ is defined as:
\[
S^{(k)} = {A^{(k)}}_{bb} - {A^{(k)}}_{bi} \, {{A^{(k)}}_{ii}}^{-1} \, {A^{(k)}}_{ib}
\]
By inverting the factorization, it immediately follows that the bottom-right block of ${A^{(k)}}^{-1}$ coincides with ${S^{(k)}}^{-1}$. Since ${A^{(k)}}^{-1}$ is the discrete analog of the solution operator of problem (\ref{pb:1}), we conclude that ${S^{(k)}}^{-1}$ is the restriction to $\Omega_b^{(k)}$ of the discrete solution operator. Consequently, the inverse Schur complement has rank-deficient off-diagonal blocks. Since such class of matrices is closed with respect to inversion, the Schur complement has rank-deficient off-diagonal blocks as well.

We employ the factorizations of the sub-problems (\ref{pb:1}) to obtain and $LDM^t$ factorization of the stiffness matrix $A$ of the original problem (\ref{pb:0}):
\[ A = 
L^{(2)}L^{(1)}
\underset{D}{
\left(
\begin{array}{cc | cc}
{A^{(1)}}_{ii}	&				& 				& 	\\
			& {A^{(2)}}_{ii}	&				& 	\\ \hline
			&				& {S^{(1)}}		& {A^{(1,2)}} 		\\
			& 				& {A^{(2,1)}}		& {S^{(2)}}
\end{array}\right)
}
{M^{(1)}}^t {M^{(2)}}^t
\]
Here $L^{(\cdot)}$ and $M^{(\cdot)}$ are accumulated Gauss transforms, and the Schur complements $S^{(\cdot)}$ are defined as previously. In order to further proceed with the factorization, we manipulate $D_{BR}$, namely the bottom right block of $D$, in the following way:
\begin{multline*}
D_{BR} = 
\left(
\begin{array}{cc}
{S^{(1)}}		& {A^{(1,2)}} 		\\
{A^{(2,1)}}	& {S^{(2)}}
\end{array}
\right)
\xrightarrow{\text{repartition}}
\left(
\begin{array}{ cc | cc}
{S^{(1)}}_{ii}		& {S^{(1)}}_{ib}	& {A^{(1,2)}}_{ii} 	& {A^{(1,2)}}_{ib} 		\\
{S^{(1)}}_{bi}		& {S^{(1)}}_{bb}	& {A^{(1,2)}}_{bi} 	& {A^{(1,2)}}_{bb} 		\\ \hline
{A^{(2,1)}}_{ii}	& {A^{(2,1)}}_{ib}	& {S^{(2)}}_{ii}	& {S^{(2)}}_{ib}		\\
 {A^{(2,1)}}_{bi}	& {A^{(2,1)}}_{bb}	& {S^{(2)}}_{bi}	& {S^{(2)}}_{bb}		\\
\end{array}\right) \\
\xrightarrow{\text{permute}}
\left(
\begin{array}{cc | cc}
{S^{(1)}}_{ii}		& {A^{(1,2)}}_{ii}	& {S^{(1)}}_{ib} 	& {A^{(1,2)}}_{ib} 		\\
{A^{(2,1)}}_{ii}	& {S^{(2)}}_{ii}	& {A^{(2,1)}}_{ib} 	& {S^{(2)}}_{ib} 		\\ \hline
{S^{(1)}}_{bi}		& {A^{(1,2)}}_{bi}	& {S^{(1)}}_{bb}	& {A^{(1,2)}}_{bb}		\\
{A^{(2,1)}}_{bi}	& {S^{(2)}}_{bi}	& {A^{(2,1)}}_{bb}	& {S^{(2)}}_{bb}
\end{array}\right)
\xrightarrow{\text{regroup}}
\left(
\begin{array}{cc}
{\hat{A}{}^{(0)}}_{ii}	& {\hat{A}^{(1,2)}} 		\\
{\hat{A}^{(2,1)}}	& {\hat{A}{}^{(0)}}_{bb}
\end{array}
\right)
\end{multline*}
Up to a permutation we shall omit, we obtain the final factorization:
\[ A = 
L^{(0)} L^{(2)} L^{(1)}
\left(
\begin{array}{cc cc}
{A^{(1)}}_{ii}	&				& 					& \\
			& {A^{(2)}}_{ii}	&					& \\
			&				& {\hat{A}{}^{(0)}}_{ii}	& \\
			& 				& 					& {S^{(0)}}
\end{array}\right)
{M^{(1)}}^t {M^{(2)}}^t {M^{(0)}}^t
\]
where the Schur complement is defined as:
\[
{S^{(0)}} = {\hat{A}{}^{(0)}}_{bb} - {\hat{A}^{(2,1)}} \, {{\hat{A}{}^{(0)}}_{ii}}^{-1} \,  {\hat{A}^{(1,2)}}
\]
If define $\overline{\Omega_b} = \cup \, \{ \supp \varphi_j \: : \: j \in  \ind  ({\hat{A}{}^{(0)}}_{bb}) \}$, the previous reasoning allows us to characterize ${S^{(0)}}^{-1}$ as the restriction of the solution operator of problem (\ref{pb:0}) to $\Omega_b$. We conclude that ${S^{(0)}}^{-1}$ and ${S^{(0)}}$ have low-rank off-diagonal blocks.

Finally, let us discuss the case of $\mathcal{L}$ being the Helmholtz operator. Although the interpretation of the inverse Schur complements as restrictions of the solution operator is preserved, the underlying Green's function not longer exhibits long-range low-rank interactions. However, in the case of scattering problems that involve relatively thin and elongated structures, the Green's function still exhibits low-rank behavior, see \cite{martinsson2007fast}. As remarked in \cite{engquist2011sweeping}, the rank is dependent on the boundary condition employed to terminate the domain, which in our scenario is $\mathcal{B} \, u = g$ on $\Gamma \cap \partial \Omega^{(k)}$, and only a PML condition guarantees consistently low ranks. Let us interpret the sets $\Omega_b^{(1)}$, $\Omega_b^{(2)}$ and $\Omega_b$ as elongated scatterers. Although we are not using a PML boundary condition, we conjecture that the observed numerical ranks are sufficiently small to yield a cheap and reliable approximation of the solution operator.

\section{Preconditioner Construction}\label{sec:3}

Our construction relies on a variant of the well-known nested dissection algorithm introduced by George in \cite{george1973nested}. We extend the work of Gilmann and Martinsson, see \cite{gillman2014n}, developed in the context of finite difference approximations of elliptic PDE's, to finite element approximations of a general differential operator. In order to do so, instead of relying on geometrical considerations similar to domain decomposition techniques, we employ a purely algebraic, black-box approach that allows us to handle a much more general framework. We remark that, although the construction of the preconditioner is fully general, we expect its performance to be affected by the nature of the underlying PDE.

Our approach rests upon a recursive reordering of the degrees of freedom (dof's) of a finite element stiffness matrix $A$. The reordering generates a binary tree which, in return, dictates a hierarchy of Schur complements. The novelty is that, provided that $A$ is sparse and that all Schur complements are to high precision Hierarchically-Semi-Separable, the following facts hold true:
\begin{enumerate}
\item an $L D M^t$ factorization of $A$ can be realized within linear complexity;
\item the factorization can be inverted within linear complexity;
\item the inverse can be applied within linear complexity.
\end{enumerate}
The factorization can be realized with a trivially parallel process, whose cost is dominated by the cost of processing the Schur complements on the top tree level. The fact that the inverse factorization can be applied fast makes it perfectly suitable to be employed as a preconditioner.

\subsection{Matrix Reordering}\label{sec:3.0}

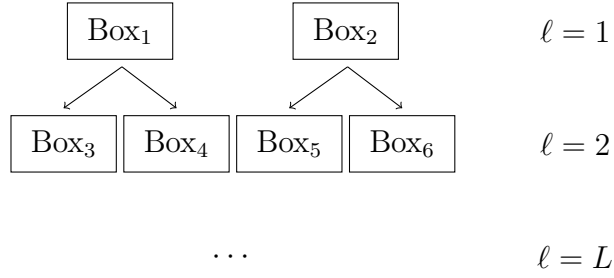
\begin{figure}[t]

\begin{center}
\begin{tikzpicture}[scale = 0.05]

\draw[black] (0,0) rectangle (28,15);
\draw[black] (30,0) rectangle (58,15);
\draw[black] (60,0) rectangle (88,15);
\draw[black] (90,0) rectangle (118,15);
\node at (14,7.5) {$\text{Box}_3$};
\node at (44,7.5) {$\text{Box}_4$};
\node at (74,7.5) {$\text{Box}_5$};
\node at (104,7.5) {$\text{Box}_6$};

\draw[black] (15,30) rectangle (43,45);
\draw[black] (75,30) rectangle (103,45);
\node at (29,37.5) {$\text{Box}_1$};
\node at (89,37.5) {$\text{Box}_2$};

\draw[black,->] (29.5,28)--(14,17);
\draw[black,->] (29.5,28)--(44,17);
\draw[black,->] (89.5,28)--(74,17);
\draw[black,->] (89.5,28)--(104,17);

\node at (150,37.5) {$\ell = 1$};
\node at (150, 7.5) {$\ell = 2$};
\node at (150, -22.5) {$\ell = L$};
\node at (59, -22.5) {$\cdots$};

\end{tikzpicture}
\end{center}
\caption{Binary tree of boxes for dof's partition}\label{fig:4}
\end{figure}

We partition the dof's of $A$ into two boxes, $\text{Box}_1, \text{Box}_2$, and, for each box $\text{Box}_i$, identify interior dof's $I^i$ and boundary dof's $B^i$ in the sense of the following connectivity graph:
\begin{equation}\label{graph:1}
\begin{tikzpicture}[<->,>=stealth',shorten >=1pt,node distance=2.5cm,auto,main node/.style={rectangle,rounded corners,draw,align=center}]

\node[main node] (1) {$B^1$};
\node[main node] (2) [right of=1] {$B^2$};

\node[main node] (3) [below of=1] {$I^1$};
\node[main node] (4) [right of=3] {$I^2$};

\path
(1) edge node [swap] {}(2)
(1) edge node [swap] {} (3)
(2) edge node [swap] {} (4);
\end{tikzpicture}
\end{equation}
The interpretation is straightforward: distinct boxes are connected to each other, in the sense of an algebraic graph, through their boundaries dof's only. At this level of exposition, the most attractive partition is the one that maximizes the number of interior dof's while minimizing the number of boundary dof's. The construction proceeds by repartitioning each box into a pair of sibling boxes, in fact creating a binary tree of boxes\footnote{Strictly speaking, this is a tree with missing route, i.e., a forest. In fact, if we were to introduce a single top-level box holding the entirety of the dof's, under a purely algebraic approach, no boundary dof's could be identified.}, see Figure \ref{fig:4}, and by identifying boundary and interior dof's, in the sense of graph (\ref{graph:1}), for each new pair of  sibling boxes. For illustration, the connectivity graph of boxes on the second tree level, i.e., $\ell = 2$, is:
\[
\begin{tikzpicture}[<->,>=stealth',shorten >=1pt,node distance=2.5cm,auto,main node/.style={rectangle,rounded corners,draw,align=center}]

\node[main node] (1) {$B^3$};
\node[main node] (2) [right of=1] {$B^4$};
\node[main node] (3) [right of=2] {$B^5$};
\node[main node] (4) [right of=3] {$B^6$};

\node[main node] (5) [below of=1] {$I^3$};
\node[main node] (6) [right of=5] {$I^4$};
\node[main node] (7) [right of=6] {$I^5$};
\node[main node] (8) [right of=7] {$I^6$};

\path
(1) edge node [swap] {}(2)
(2) edge node [swap] {} (3)
(3) edge node [swap] {} (4)
(1) edge node [swap] {} (5)
(2) edge node [swap] {} (6)
(3) edge node [swap] {} (7)
(4) edge node [swap] {} (8)
(1) edge [bend left=40] node  {} (3)
(1) edge [bend left=40] node  {} (4)
(2) edge [bend left=40] node  {} (4);

\end{tikzpicture}
\]
The construction terminates at tree level $\ell = L$, when the newly created boxes contain a number of dof's sufficiently small to allow for dense linear algebra operations at a negligible cost. For consistency with the notation introduced in Section \ref{sec:1}, we switch to Greek letters and identify a box with its index, i.e., $\sigma = \text{Box}_\sigma$.

We establish an ordering for the dof's in $I^\sigma$ and $B^\sigma$, and define the following sub-matrices of $A$:
\begin{align*}
{A^{(\sigma)}}_{ii}	& = A(I^\sigma,I^\sigma) 		&	\qquad 	& \sigma\text{-box}\text{ interior-to-interior}  \\
{A^{(\sigma)}}_{bb}	& = A(B^\sigma,B^\sigma)  	&			& \sigma\text{-box}\text{ boundary-to-boundary}  \\
{A^{(\sigma)}}_{bi} 	& = A(B^\sigma,I^\sigma) 	&			& \sigma\text{-box}\text{ boundary-to-interior} \\
{A^{(\sigma)}}_{ib} 	& = A(I^\sigma,B^\sigma) 	&			& \sigma\text{-box}\text{ interior-to-boundary}\\
{A^{(\sigma,\tau)}} 	& = A(B^\sigma,B^\tau) 		&			& \sigma\text{-box}\text{ to } \tau \text{-box} \text{ (boundary interaction only)}
\end{align*}
Let $\sigma_1 , \ldots , \sigma_n$ be the leaf boxes and order the dof's of $A$ by grouping together the interior dof's $I^{\sigma_1}, \ldots , I^{\sigma_n}$, and the boundary dof's $B^{\sigma_1} , \ldots , B^{\sigma_n}$ on tree level $\ell = L$:
\begin{equation}\label{eq:30}
A = 
\begin{pmatrix}[ccc | ccc]
{A^{(\sigma_1)}}_{ii}	&				&					& {A^{(\sigma_1)}}_{ib}	&				& 				\\
					& \ddots 			&					& 					& \ddots			& 				\\
					&				& {A^{(\sigma_n)}}_{ii}	&					&				& {A^{(\sigma_n)}}_{ib}	\\ \hline
{A^{(\sigma_1)}}_{bi}	&				&					& {A^{(\sigma_1)}}_{bb}	& \cdots			& {A^{(\sigma_1,\sigma_n)}} 		\\
					& \ddots			&					& \vdots				& \ddots			& \vdots			\\
					&				& {A^{(\sigma_n)}}_{bi}	& {A^{(\sigma_n,\sigma_1)}}	& \cdots			&  {A^{(\sigma_n)}}_{bb}
\end{pmatrix}
\end{equation}
As common practice, we have omitted the null blocks.

\subsection{Hierarchical Matrix Factorization}\label{sec:3.1}

The factorization strategy recursively decouples interior dof's from boundary dof's through Gauss transforms, starting from the leaf-boxes, all the way up to the top tree-level. The Schur complements that arise in the process are treated through accelerated linear algebra techniques. More specifically, the complement associated to a box is obtained by a fast merge of the complements of its child-boxes.

For each box $\sigma$, we define Gauss transforms as the unit-lower-triangular matrices $L^{(\sigma)}$ and $M^{(\sigma)}$ such that:
\begin{equation}
L^{(\sigma)}(B^\sigma, I^\sigma) = {A^{(\sigma)}}_{bi} \,{{A^{(\sigma)}}_{ii}}^{-1} \qquad , \qquad M^{(\sigma)}(B^\sigma, I^\sigma) = \Big( {{A^{(\sigma)}}_{ii}}^{-1} \, {A^{(\sigma)}}_{ib}\Big)^{t}
\end{equation}
Since it is well-understood how Gauss transforms accumulate and commute with permutations, we henceforth refer to them generically as $L$ and $M$. We employ the Gauss transforms to decouple the top-left super-block of $A$:
\begin{equation}\label{eq:6}
L^{-1} A M^{-t} = 
\begin{pmatrix}[ ccc|c cccc ]
{A^{({\sigma}_1)}}_{ii}	&			&						& 					&					&					& 					\\
					& \ddots		&						& 					& 					&					& 					\\
					&			& {A^{({\sigma}_n)}}_{ii}	&					&					&					&  					\\ \hline
					&			&						& {S^{({\sigma}_1)}}		& {A^{({\sigma}_1,{\sigma}_2)}}	& \cdots				& {A^{({\sigma}_1,{\sigma}_n)}} 	\\
					& 			&						& {A^{({\sigma}_2,{\sigma}_1)}}	& \ddots				& \ddots				& \vdots				\\
					&			&						& \vdots				& \ddots				& \ddots				& {A^{({\sigma}_{n-1},{\sigma}_n)}} 	\\
					&			& 						& {A^{({\sigma}_n,{\sigma}_1)}}	& \cdots				& {A^{({\sigma}_n,{\sigma}_{n-1})}}	& {S^{({\sigma}_n)}}
\end{pmatrix}
\end{equation}
The top-right and bottom-left super-blocks vanish, while Schur complements $S^{(\sigma)} = {A^{(\sigma)}}_{bb} -  {A^{(\sigma)}}_{bi} \, {{A^{(\sigma)}}_{ii}}^{-1} \, {A^{(\sigma)}}_{ib}$ appear on the diagonal of the bottom-right super-block. 

In order to recursively proceed in the factorization, for each pair of sibling boxes $\{ \mu, \nu \}$ with common ancestor $\sigma$, we define the remaining interior dof's $\hat{I}$ and the remaining boundary dof's $\hat{B}^\sigma$ as:
\[
\hat{I}^\sigma = (B^\mu \cup B^\nu) \cap I^\sigma \qquad , \qquad \hat{B}^\sigma =   (B^\mu \cup B^\nu) \cap B^\sigma
\]
Consequently, $\{\hat{I}^\sigma , \hat{B}^\sigma \}$ is a partition of the aggregated boundary $B^\mu \cup B^\nu$, while $\{ B^\mu \cap \hat{I}^\sigma , B^\mu \cap \hat{B}^\sigma \}$ is a partition\footnote{The two sets are evidently disjoint. Furthermore: $(B^\mu \cap \hat{I}^\sigma) \cup (B^\mu \cap \hat{B}^\sigma) = B^\mu \cap ( \hat{I}^\sigma \cup \hat{B}^\sigma) =  B^\mu \cap ( B^\mu \cup B^\nu) = B^\mu$.} of $B^\mu$. Up to a permutation we shall omit, we partition the Schur complement $S^{(\mu)}$ as:
\[
S^{(\mu)} = 
\begin{pmatrix}
{S^{(\mu)}}_{ii}	& {S^{(\mu)}}_{ib}	\\
{S^{(\mu)}}_{bi}	& {S^{(\mu)}}_{bb}
\end{pmatrix}
\]
where the blocks are defined as:
\begin{align*}
{S^{(\mu)}}_{ii}	& = S^{(\mu)}(B^\mu \cap \hat{I}^\sigma,B^\mu \cap \hat{I}^\sigma)	\\
{S^{(\mu)}}_{bb}	& = S^{(\mu)}(B^\mu \cap \hat{B}^\sigma,B^\mu \cap \hat{B}^\sigma)	\\
{S^{(\mu)}}_{bi} 	& = S^{(\mu)}(B^\mu \cap \hat{B}^\sigma,B^\mu \cap \hat{I}^\sigma)	\\
{S^{(\mu)}}_{ib} 	& = S^{(\mu)}(B^\mu \cap \hat{I}^\sigma,B^\mu \cap \hat{B}^\sigma)
\end{align*}
Let us define matrices:
\begin{equation*}
{\hat{A}{}^{(\sigma)}}_{ii} = 
\begin{pmatrix}
{S^{(\mu)}}_{ii} & \hat{A}{}^{(\mu,\nu)} \\
\hat{A}{}^{(\nu,\mu)} & S^{{(\nu)}_{ii}}
\end{pmatrix}
\: ; \:
{\hat{A}{}^{(\sigma)}}_{ib} = 
\begin{pmatrix}
{S^{(\mu)}}_{ib} & \hat{A}{}^{(\mu,\nu)} \\
\hat{A}{}^{(\nu,\mu)} & {S^{(\nu)}}_{ib}
\end{pmatrix}
\: ; \:
{\hat{A}{}^{(\sigma)}}_{bb} = 
\begin{pmatrix}
{S^{(\mu)}}_{bb} & \hat{A}{}^{(\mu,\nu)} \\
\hat{A}{}^{(\nu,\mu)} & {S^{(\nu)}}_{bb}
\end{pmatrix}
\end{equation*}
and ${\hat{A}{}^{(\sigma)}}_{bi}$ analogously to ${\hat{A}{}^{(\sigma)}}_{ib}$. In order to simplify the notation, $\hat{A}^{(\mu,\nu)}$ indicates a generic sub-matrix of $A^{(\mu,\nu)}$ that can be inferred from the context\footnote{Distinct instances of the symbol should be regarded as different matrices.}. By  reordering the boundary dof's $B^{\sigma_1} , \ldots , B^{\sigma_n}$ as $\hat{I}^{(\cdot)} , \ldots , \hat{I}^{(\cdot)} , \hat{B}^{(\cdot)} , \ldots , \hat{B}^{(\cdot)}$, where the omitted superscripts are the parent boxes of the leaves $\sigma_1, \ldots , \sigma_n$, equation (\ref{eq:6}) becomes:
\[
L^{-1} A M^{-t} = 
\begin{pmatrix}[c | ccc  ccc]
\star & & & & & & \\ \hline
& {\hat{A}{}^{(\cdot)}}_{ii}	&				&						& {\hat{A}{}^{(\cdot)}}_{ib}	&				& 						\\
&						& \ddots 			&						& 						& \ddots			& 						\\
&						&				& {\hat{A}{}^{(\cdot)}}_{ii}	&						&				&{\hat{A}{}^{(\cdot)}}_{ib}	\\
&{\hat{A}{}^{(\cdot)}}_{bi}	&				&						& {\hat{A}{}^{(\cdot)}}_{bb}	& \cdots			& \hat{A}{}^{(\cdot,\cdot)} 	\\
&						& \ddots			&						& \vdots					& \ddots			& \vdots					\\
&						&				& {\hat{A}{}^{(\cdot)}}_{bi}	& \hat{A}{}^{(\cdot,\cdot)}	& \cdots			&  {\hat{A}{}^{(\cdot)}}_{bb}
\end{pmatrix}
\]
where the symbol $\star$ indicates an omitted non-zero block.

The factorization proceeds recursively by exploiting the fact that the bottom-right super-block has structure that is identical to that of the original matrix $A$, compare to equation (\ref{eq:30}). The interior dof's $\hat{I}^{(\cdot)} , \ldots , \hat{I}^{(\cdot)}$ are recursively decoupled through Gauss transforms until the top tree-level is reached. The Schur complement that originates for the elimination of the $\hat{I}^\sigma$ dof's is defined as:
\begin{equation}\label{eq:32}
S^{(\sigma)}
 = 
\underset{{\hat{A}{}^{(\sigma)}}_{bb}}{
\begin{pmatrix}
{S^{(\mu)}}_{bb} 		& {\hat{A}{}^{(\mu,\nu)}} \\
{\hat{A}{}^{(\nu,\mu)}}	& {S^{(\nu)}}_{bb} 	
\end{pmatrix}
} -
\underset{{\hat{A}{}^{(\sigma)}}_{bi}}{
\begin{pmatrix}
{S^{(\mu)}}_{bi} 		& {\hat{A}{}^{(\mu,\nu)}} \\
{\hat{A}{}^{(\nu,\mu)}}	& {S^{(\nu)}}_{bi} 	
\end{pmatrix}
}
\underset{{{\hat{A}{}^{(\sigma)}}_{ii}}}{
\begin{pmatrix}
{S^{(\mu)}}_{ii} 		& {\hat{A}{}^{(\mu,\nu)}} \\
{\hat{A}{}^{(\nu,\mu)}}	& {S^{(\nu)}}_{ii} 	
\end{pmatrix}
}^{-1}
\underset{{\hat{A}{}^{(\sigma)}}_{ib}}{
\begin{pmatrix}
{S^{(\mu)}}_{ib} 		& {\hat{A}{}^{(\mu,\nu)}} \\
{\hat{A}{}^{(\nu,\mu)}}	& {S^{(\nu)}}_{ib} 	
\end{pmatrix}
}
\end{equation}
When the decoupling process terminates, we obtain the following factorization:
\begin{equation}\label{eq:31}
L^{-1} A  M^{-t} = 
\begin{pmatrix}[ ccc cc ]
{\hat{A}{}^{(\sigma_1)}}_{ii}	& 			 	& 							&				&			\\
 						& \ddots 			&							&				&			\\ 
						& 				& {\hat{A}{}^{(\sigma_n)}}_{ii}	& 				&			\\
						&				&							& S^{(1)}			& A^{(1,2)} 	\\
						&				&							& A^{(2,1)}		& S^{(2)}
\end{pmatrix}
\end{equation}
The boxes $\sigma_1, \ldots, \sigma_n$ are ordered starting from the bottom tree-level and, for consistency of notation, we have set ${\hat{A}{}^{(\sigma)}}_{ii} = {A^{(\sigma)}}_{ii}$ for the leaf-boxes as well. The matrices $L$ and $M$ are the accumulated Gauss transforms\footnote{As before, permutations have been omitted in the definition of the Gauss transforms.}, relative to all boxes, excluding those on the top tree-level:
\begin{equation*}
L = L^{(\sigma_1)} \cdots L^{(\sigma_n)} \qquad ; \qquad M  = M^{(\sigma_1)} \cdots M^{(\sigma_n)}
\end{equation*}
By setting
\begin{equation}\label{eq:36}
\hat{A}^{(0)} = 
\begin{pmatrix}
S^{(1)}	& A^{(1,2)} \\
A^{(2,1)}	& S^{(2)}
\end{pmatrix} \qquad ; \qquad
D = \diag \{ {\hat{A}{}^{(\sigma_1)}}_{ii} , \ldots , {\hat{A}{}^{(\sigma_n)}}_{ii} , \hat{A}^{(0)} \}
\end{equation}
we obtain the desired $LDM^t$ factorization of $A$.

\subsection{Fast Merging of Schur Complements}
The hierarchy of the boxes dictates a hierarchy of the Schur complements in the sense that $S^{(\sigma)}$ depends only upon the Schur complements of the child-boxes of $\sigma$, namely $S^{(\mu)}$ and $S^{(\nu)}$,  and some blocks of the original matrix $A$, see equation (\ref{eq:32}). Informally, we say that $S^{(\sigma)}$ is obtained by ``merging'' together $S^{(\mu)}$ and $S^{(\nu)}$. As we are about to show, when $S^{(\mu)}$ and $S^{(\nu)}$ are in HSS-form, the merging procedure can be performed fast, in the sense that an HSS approximation to $S^{(\sigma)}$ can be computed within linear complexity.

The fast merging procedure is based the following assumptions:
\begin{enumerate}
\item all sub-matrices $\hat{A}^{(\cdot , \cdot)}$ of $A$ are sparse; \label{en:item:1}
\item all leaf-node Schur complements $S^{(\cdot )}$ and their inverses ${S^{(\cdot )}}^{-1}$ are in HSS-form.
\end{enumerate}
In general, the property of a matrix to be sparse does not carry over to an arbitrarily selected sub-matrix. In the context of high-order finite element approximations, this is well-established. In fact, dense blocks allow advanced finite element solvers to employ high-performance dense linear algebra. We remark that assumption (\ref{en:item:1}) is much milder, since it only requires matrices describing boundary-to-boundary sub-interactions to be sparse.

For ease of exposition, let us recall the definition of $S^{(\sigma)}$ as in (\ref{eq:32}):
\[
S^{(\sigma)}
 = 
\underset{{\hat{A}{}^{(\sigma)}}_{bb}}{
\begin{pmatrix}
{S^{(\mu)}}_{bb} 		& {\hat{A}{}^{(\mu,\nu)}} \\
{\hat{A}{}^{(\nu,\mu)}}	& {S^{(\nu)}}_{bb} 	
\end{pmatrix}
} -
\underset{{\hat{A}{}^{(\sigma)}}_{bi}}{
\begin{pmatrix}
{S^{(\mu)}}_{bi} 		& {\hat{A}{}^{(\mu,\nu)}} \\
{\hat{A}{}^{(\nu,\mu)}}	& {S^{(\nu)}}_{bi} 	
\end{pmatrix}
}
\underset{{{\hat{A}{}^{(\sigma)}}_{ii}}}{
\begin{pmatrix}
{S^{(\mu)}}_{ii} 		& {\hat{A}{}^{(\mu,\nu)}} \\
{\hat{A}{}^{(\nu,\mu)}}	& {S^{(\nu)}}_{ii} 	
\end{pmatrix}
}^{-1}
\underset{{\hat{A}{}^{(\sigma)}}_{ib}}{
\begin{pmatrix}
{S^{(\mu)}}_{ib} 		& {\hat{A}{}^{(\mu,\nu)}} \\
{\hat{A}{}^{(\nu,\mu)}}	& {S^{(\nu)}}_{ib} 	
\end{pmatrix}
}
\]
Assume that $S^{(\mu)}$ and $S^{(\nu)}$ are in HSS-form. The fast merging procedure is performed through the following steps.
\begin{description}[style=multiline,leftmargin=2cm,]
\item[Step 1]
Since a sub-matrix of an HSS matrix can be trivially obtained through a tall and skinny permutation matrix $P$, e.g., ${S^{(\mu)}}_{bb} = P^t \, S^{(\mu)} \, P$, matrices ${\hat{A}{}^{(\sigma)}}_{bb}$, ${\hat{A}{}^{(\sigma)}}_{bi}$, and ${\hat{A}{}^{(\sigma)}}_{ib}$ can be applied to a vector within linear complexity.
\item[Step 2]
The action of the inverse of ${\hat{A}{}^{(\sigma)}}_{ii}$ on a vector $z$ is 
equivalent to the solution of the linear system:
\[
\begin{pmatrix}
{S^{(\mu)}}_{ii} 		& \hat{A}^{(\mu,\nu)} \\
\hat{A}^{(\nu,\mu)}	& {S^{(\nu)}}_{ii} 	
\end{pmatrix}
\begin{pmatrix}
x_1\\
x_2
\end{pmatrix}=
\begin{pmatrix}
z_1\\
z_2
\end{pmatrix}
\]
where $z = (z_1,z_2)$ and $x = (x_1,x_2)$ have been partitioned according to the blocking of ${\hat{A}{}^{(\sigma)}}_{ii}$. A standard block-solve yields:
\begin{subequations}\label{eq:37}
\begin{align}
x_2 & = {{\tilde{S}^{(\nu)}{}}_{ii}}^{-1} (z_2 - \hat{A}^{(\nu,\mu)} \, {{S^{(\mu)}}_{ii}}^{-1} \, z_1) \\
x_1 & = {{S^{(\mu)}}_{ii}}^{-1} \, z_1 - {{S^{(\mu)}}_{ii}}^{-1} \, \hat{A}^{(\mu,\nu)} \, x_2
\end{align}
\end{subequations}
where ${{\tilde{S}^{(\nu)}{}}_{ii}} = {S^{(\nu)}}_{ii} - \hat{A}^{(\nu,\mu)} \, {{S^{(\mu)}}_{ii}}^{-1} \, \hat{A}^{(\mu,\nu)}$. Let us remark that matrices ${{S^{(\mu)}}_{ii}}^{-1}$ and $ {{\tilde{S}^{(\nu)}{}}_{ii}}^{-1}$ fully describe the action of ${{\hat{A}{}^{(\sigma)}}_{ii}}^{-1}$ through equations (\ref{eq:37}). Since ${{\tilde{S}^{(\nu)}{}}_{ii}}$ can be applied within linear complexity, by virtue of Theorem \ref{th:1}, it can be compressed efficiently. The cost of the current step is:
\[
\mathrm{O}(\# \hat{I}^\sigma k^2) \quad = \quad
\underset{\text{inversion of  ${S^{(\mu)}}_{ii}$}}{\mathrm{O}(\# \hat{I}^\sigma k^2)} +
\underset{\text{compression of  ${{\tilde{S}^{(\nu)}{}}_{ii}}$}}{\mathrm{O}(\# \hat{I}^\sigma k^2)} +
\underset{\text{inversion of  ${{\tilde{S}^{(\nu)}{}}_{ii}}$}}{\mathrm{O}(\# \hat{I}^\sigma k^2)} 
\]
The procedure is detailed in Algorithm (\ref{algo:1}). 
\item[Step 3]
Under the assumption that $\# \hat{B}^\sigma = \mathrm{O}(\# \hat{I}^\sigma)$, the previous steps imply that $S^{(\sigma)}$ can be applied to a vector within linear complexity. By Theorem \ref{th:1}, it can be reduced to HSS form at the cost $\mathrm{O}(\# \hat{B}^\sigma k^2)$.
\end{description}

\begin{algorithm}[t]
  \SetAlgoLined
\SetKwInOut{Input}{input}\SetKwInOut{Output}{output}
  \Input{$z_1$, $z_2$, $S^{(\mu)}$, $S^{(\nu)}$, $\hat{A}^{(\mu,\nu)}$, $\hat{A}^{(\nu,\mu)}$}
  \Output{$x_1$, $x_2$, ${{S^{(\mu)}}_{ii}}^{-1}$, ${{\tilde{S}^{(\nu)}{}}_{ii}}^{-1}$}
\vspace{0.2cm}
\texttt{-- STEP 2.1 --}\newline
determine permutations $P_{\mu} : B^{\mu} \to B^{\mu} \cap \hat{I}^{\sigma}$, and $P_{\nu} : B^{\nu}\to B^{\nu} \cap \hat{I}^{\sigma}$\;
define ${S^{(\mu)}}_{ii} = P^t_{\mu} \, S^{(\mu)} \, P_{\mu}$, and ${S^{(\nu)}}_{ii} = P^t_{\nu} \, S^{(\nu)} \, P_{\nu}$\;
compute ${{S^{(\mu)}}_{ii}}^{-1}$ through fast inversion\;
\vspace{0.2cm}
\texttt{-- STEP 2.2 --}\newline
compress ${{\tilde{S}^{(\nu)}{}}_{ii}} = {S^{(\nu)}}_{ii} - \hat{A}^{(\nu,\mu)} \, {{S^{(\mu)}}_{ii}}^{-1} \, \hat{A}^{(\mu,\nu)}$ as in Theorem \ref{th:1}\;
compute ${{\tilde{S}^{(\nu)}{}}_{ii}}^{-1}$ through fast inversion\;
\vspace{0.2cm}
\texttt{-- STEP 2.3 --}\newline
compute $x_2 = {{\tilde{S}^{(\nu)}{}}_{ii}}^{-1} (z_2 - \hat{A}^{(\nu,\mu)} \, {{S^{(\mu)}}_{ii}}^{-1} \, z_1)$, and $x_1  = {{S^{(\mu)}}_{ii}}^{-1} \, z_1 - {{S^{(\mu)}}_{ii}}^{-1} \, \hat{A}^{(\mu,\nu)} \, x_2$ through fast application\;
  \caption{Compression and fast application of ${{\hat{A}{}^{(\sigma)}}_{ii}}^{-1}$.}\label{algo:1}
\end{algorithm}

Let us remark that, although computing the HSS form of a parent Schur complement is a linear complexity operation, there is not guarantee that it is HSS to high precision. We expect this property to be dictated by the nature of the underlying PDE and the discretization method, and affect the cost/effectiveness ratio of the preconditioner. This is investigated through numerical experiments, that are presented in Section \ref{sec:4}.

\subsection{Approximate Matrix Factorization}\label{sec:3.2}
Let $A$ have dimension $N$ and select a number of tree levels $L$ so that each leaf box holds a sufficiently small number of dof's $m = N/2^L$ to allow for dense linear algebra manipulations at negligible cost. At the level of leaf-boxes, the Schur complements are computed and compressed to HSS-form in a straightforward manner at a cost proportional to $\mathrm{O}(m^3)$. Starting from the leaf-boxes, the Schur complements are merged together using the strategy described above, until the top tree-level is reached. The procedure is detailed in Algorithm \ref{algo:2}. The cost of the algorithm is dominated by the cost of processing the boxes on the top level, namely $\mathrm{O}(\#\hat{B}^\sigma k^2)$, $\sigma$ on level $\ell = 1$. Apart from the leaf nodes, the uncompressed Schur complements are never formed explicitly. Since the Gauss transforms $L$ and $M$ are obtained from the ${{\hat{A}{}^{(\sigma)}}_{ii}}^{-1}$ matrices, the fast merging process described in Algorithm \ref{algo:2} does, in fact, produce an approximate $LDM^t$ factorization of $A$.

\begin{algorithm}[t]
  \SetAlgoLined
\SetKwInOut{Input}{input}\SetKwInOut{Output}{output}
  \Input{$A$, binary tree partitioning the dof's of $A$}
  \Output{${{\hat{A}{}^{(\sigma)}}_{ii}}^{-1}$, $S^{(\sigma)}$ for all $\sigma$'s in the binary tree}
\vspace{0.2cm}
\texttt{-- STEP 1 --}\newline
\For{$\sigma$ on level $L$}
{
compute $S^{(\sigma)}$ by a straightforward approach\;
compress $S^{(\sigma)}$ by a straightforward approach\;
}
\vspace{0.2cm}
\texttt{-- STEP 2 --}\newline
\For{$\ell = L-1,\ldots, 1$}{
	\For{$\sigma$ on level $\ell$}{
		form ${\hat{A}{}^{(\sigma)}}_{bb}$, ${\hat{A}{}^{(\sigma)}}_{bi}$, ${\hat{A}{}^{(\sigma)}}_{ii}$, ${\hat{A}{}^{(\sigma)}}_{ib}$\;
		compress ${{\hat{A}{}^{(\sigma)}}_{ii}}^{-1}$ as in Algorithm \ref{algo:1}\;
		compress $S^{(\sigma)} = {\hat{A}{}^{(\sigma)}}_{bb} - {\hat{A}{}^{(\sigma)}}_{bi} \, {{\hat{A}{}^{(\sigma)}}_{ii}}^{-1} \, {\hat{A}{}^{(\sigma)}}_{ib}$ as in Theorem \ref{th:1}\;
	}
}
\caption{Computation of Schur complements $S^{(\sigma)}$ through fast merging.}\label{algo:2}
\end{algorithm}

Since the Gauss transforms can be trivially inverted, the cost of inverting the $LDM^t$ factorization is tantamount to the cost of inverting $D$ which, in turn, see equation (\ref{eq:36}), coincides with the cost of inverting block $\hat{A}^{(0)}$, defined as:
\[
\hat{A}^{(0)} = 
\begin{pmatrix}
S^{(1)}	& A^{(1,2)} \\
A^{(2,1)}	& S^{(2)}
\end{pmatrix}
\]
The inversion of $\hat{A}^{(0)}$ can be achieved with a slight modification of Algorithm \ref{algo:1}, described in Algorithm \ref{algo:3}. We conclude that the cost\footnote{As previously, we assume $\# \hat{B}^\sigma = \mathrm{O}(\# \hat{I}^\sigma)$.} of building the inverse factorization is
\[
\mathrm{O}(\#\hat{B}^\sigma k^2) \quad = \quad
\underset{\text{$LDM^t$ factorization}}{\mathrm{O}(\#\hat{B}^\sigma k^2)} +
\underset{\text{inversion of $\hat{A}^{(0)}$}}{\mathrm{O}(\#\hat{B}^\sigma k^2)}
 \]
 for $\sigma$ on level $\ell = 1$.

\begin{algorithm}[t]
  \SetAlgoLined
\SetKwInOut{Input}{input}\SetKwInOut{Output}{output}
  \Input{$z_1$, $z_2$, $S^{(1)}$, $S^{(2)}$, $\hat{A}^{(1,2)}$, $\hat{A}^{(2,1)}$}
  \Output{$x_1$, $x_2$, ${{S^{(1)}}}^{-1}$, ${{\tilde{S}^{(2)}{}}}^{-1}$}
\vspace{0.2cm}
\texttt{-- STEP 1 --}\newline
compute ${{S^{(1)}}}^{-1}$ through fast inversion\;
compress ${{\tilde{S}^{(2)}{}}} = {S^{(2)}} - \hat{A}^{(2,1)} \, {{S^{(1)}}}^{-1} \, \hat{A}^{(1,2)}$ as in Theorem \ref{th:1}\;
compute ${{\tilde{S}^{(2)}{}}_{ii}}^{-1}$ through fast inversion\;
\vspace{0.2cm}
\texttt{-- STEP 2 --}\newline
compute $x_2 = {{\tilde{S}^{(2)}{}}}^{-1} (z_2 - \hat{A}^{(2,1)} \, {{S^{(1)}}}^{-1} \, z_1)$, and $x_1  = {{S^{(1)}}}^{-1} \, z_1 - {{S^{(1)}}}^{-1} \, \hat{A}^{(1,2)} \, x_2$ through fast application\;
  \caption{Compression and fast application of ${\hat{A}^{(0)}{}}^{-1}$.}\label{algo:3}
\end{algorithm}

In order to obtain cost estimates with respect to the problem size $N$, we evaluate the number of dof's in $B^\sigma$ through a geometrical argument, namely:
\begin{equation}\label{eq:42}
\# B^{\sigma} = \bigg(\frac{N}{2^\ell}\bigg)^{1 - 1/d}\quad , \quad \text{for box } \sigma \text{ on level }\ell
\end{equation}
for some parameter $d$. The most conservative estimate is obtained for the limit case $d \uparrow \infty$, which yields:
\[
\# B^{\sigma} = \frac{N}{2^\ell}\quad , \quad \text{for box } \sigma \text{ on level }\ell
\]
Taking into account that $\# \hat{B}^\sigma \le \# B^\sigma$, the cost of building the inverse factorization is
\begin{equation}\label{eq:38}
\mathrm{O}(\# \hat{B}^\sigma k^2) = \mathrm{O}(N k^2)
\end{equation}
We conclude that the construction cost of the preconditioner, i.e., the inverse factorization, is linear with respect to $N$, provided that $k$ is independent of $N$.

Let us briefly comment on the last statement. The assumption that $k$ is independent of $N$, and solely dictated by the nature of the problem and the discretization, is unrealistic. In practice, we are interested in the asymptotic behavior of $k$, and some mild dependency upon $N$ is acceptable, as long as $k/N \downarrow 0$. In this respect, estimate (\ref{eq:38}) is deceptively optimistic. On the other hand, taking $d \uparrow \infty$ in estimate (\ref{eq:42}) could be excessively conservative, and lead to an overly pessimistic result. As numerical experiments presented in Section \ref{sec:4} suggest, both of those scenarios are possible and the construction cost is tightly related to strategy used to increase the problem size. We defer further discussion to Section \ref{sec:4}.

Finally, we remark that the inverse factorization $M^{-t} D^{-1} L^{-1}$ can be applied to a vector within linear complexity, thus making it suitable to be employed as a preconditioner.

\section{Numerical Results}\label{sec:4}

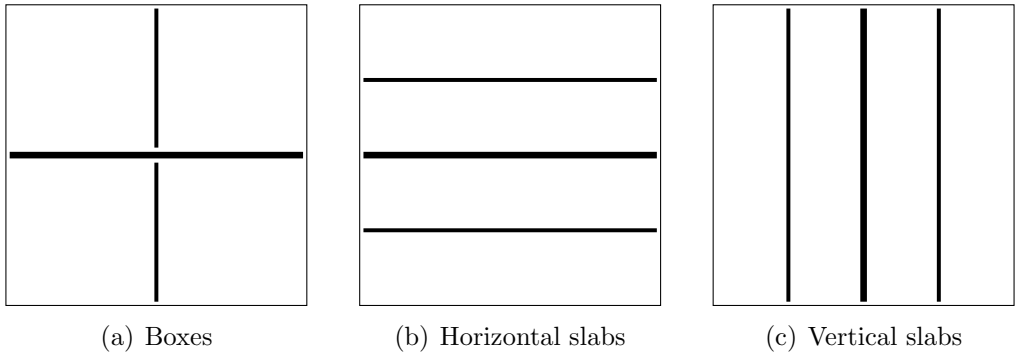
\begin{figure}[t]

\begin{center}

\subfigure[Boxes]{
\begin{tikzpicture}[scale = 0.05]

\draw[black] (0,0) rectangle (80,80);
\draw[line width = 2.5](1,40) -- (79,40);

\draw[line width = 1.5](40,1) -- (40,38);
\draw[line width = 1.5](40,42) -- (40,79);

\end{tikzpicture}\label{fig:5.1}
}\quad
\subfigure[Horizontal slabs]{
\begin{tikzpicture}[scale = 0.05]

\draw[black] (0,0) rectangle (80,80);
\draw[line width = 2.5](1,40) -- (79,40);

\draw[line width = 1.5](1,20) -- (79,20);
\draw[line width = 1.5](1,60) -- (79,60);

\end{tikzpicture}\label{fig:5.2}
}\quad
\subfigure[Vertical slabs]{
\begin{tikzpicture}[scale = 0.05]

\draw[black] (0,0) rectangle (80,80);

\draw[line width = 2.5](40,1) -- (40,79);

\draw[line width = 1.5](20,1) -- (20,79);
\draw[line width = 1.5](60,1) -- (60,79);

\end{tikzpicture}\label{fig:5.3}
}

\end{center}
\caption{Domain partitioning schemes.}\label{fig:40}
\end{figure}

Discontinuous Galerkin (DG) finite element approximations provide a natural environment for testing the proposed low-rank preconditioner. In fact, due to the doubling of dof's on the elements boundaries, the partitioning can be interpreted in purely geometrical terms, and the identification of interior and boundary dof's for each box is easily achieved. We investigate the behavior of the preconditioner for the following problems:
\begin{enumerate}
\item Poisson equation;
\item anisotropic reaction-diffusion equation;
\item Helmholtz equation;
\item Helmholtz equation with high material contrast.
\end{enumerate}
We rely upon Interior Penalty DG formulations and employ nodal DG finite element discretizations, see \cite{hesthaven2007nodal}. All problems are formulated on the unit square $\Omega = (0,1)^2$, on which we lay a structured triangular grid. We shall investigate the behavior of the preconditioner with respect to the domain partitioning scheme, see Figure \ref{fig:40}, the mesh size $h$, and the polynomial degree of approximation $p$. For the Helmholtz operator, the wave number $\kappa$ is chosen as a function of the discretization parameters $h$ and $p$, so that the dispersion error, see \cite{ainsworth2004dispersive, ainsworth2006dispersive}, remains constant as the size of the problem varies.

\begin{figure}[t]
\begin{center}
\includegraphics[scale=0.35]{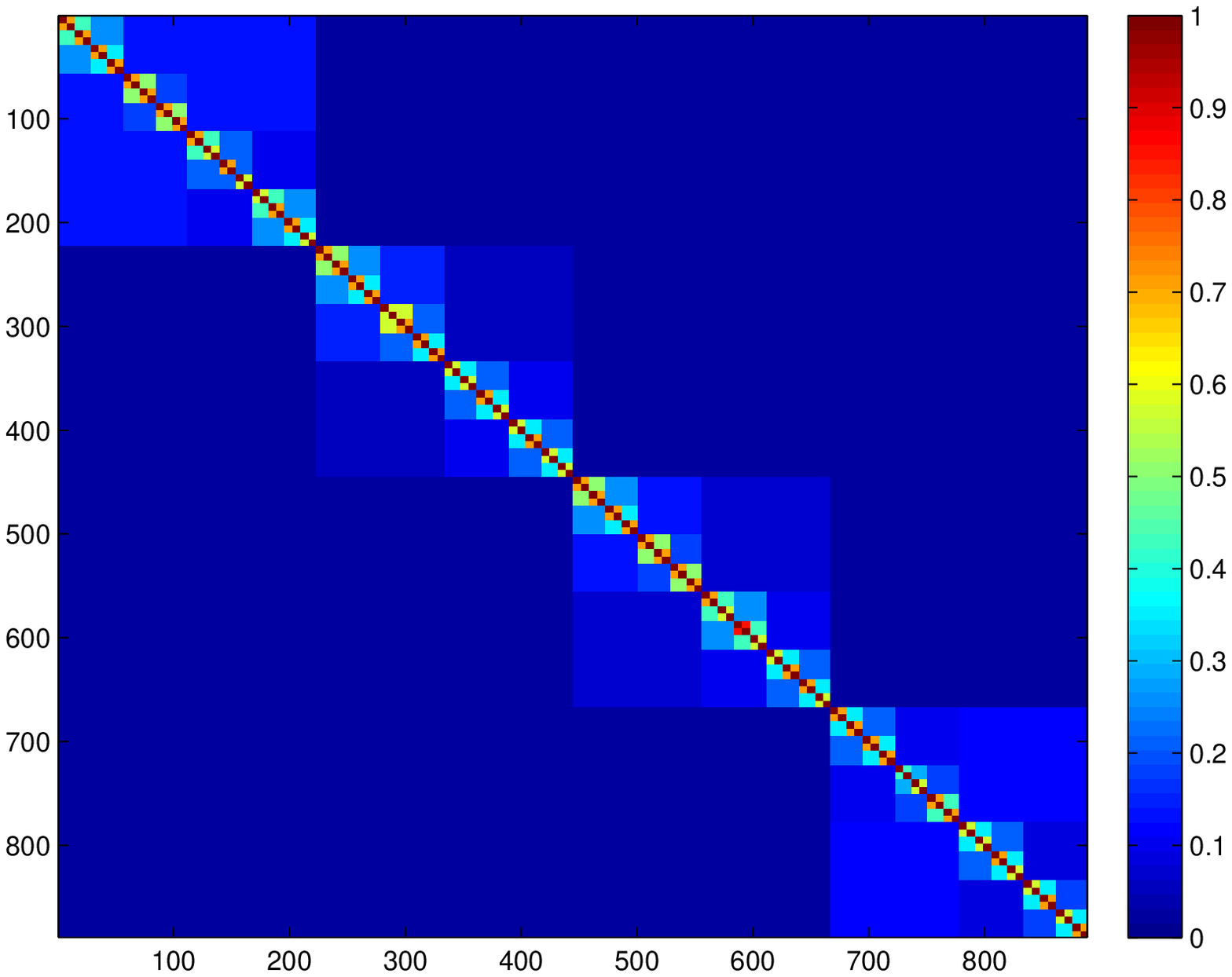}
\includegraphics[scale=0.35]{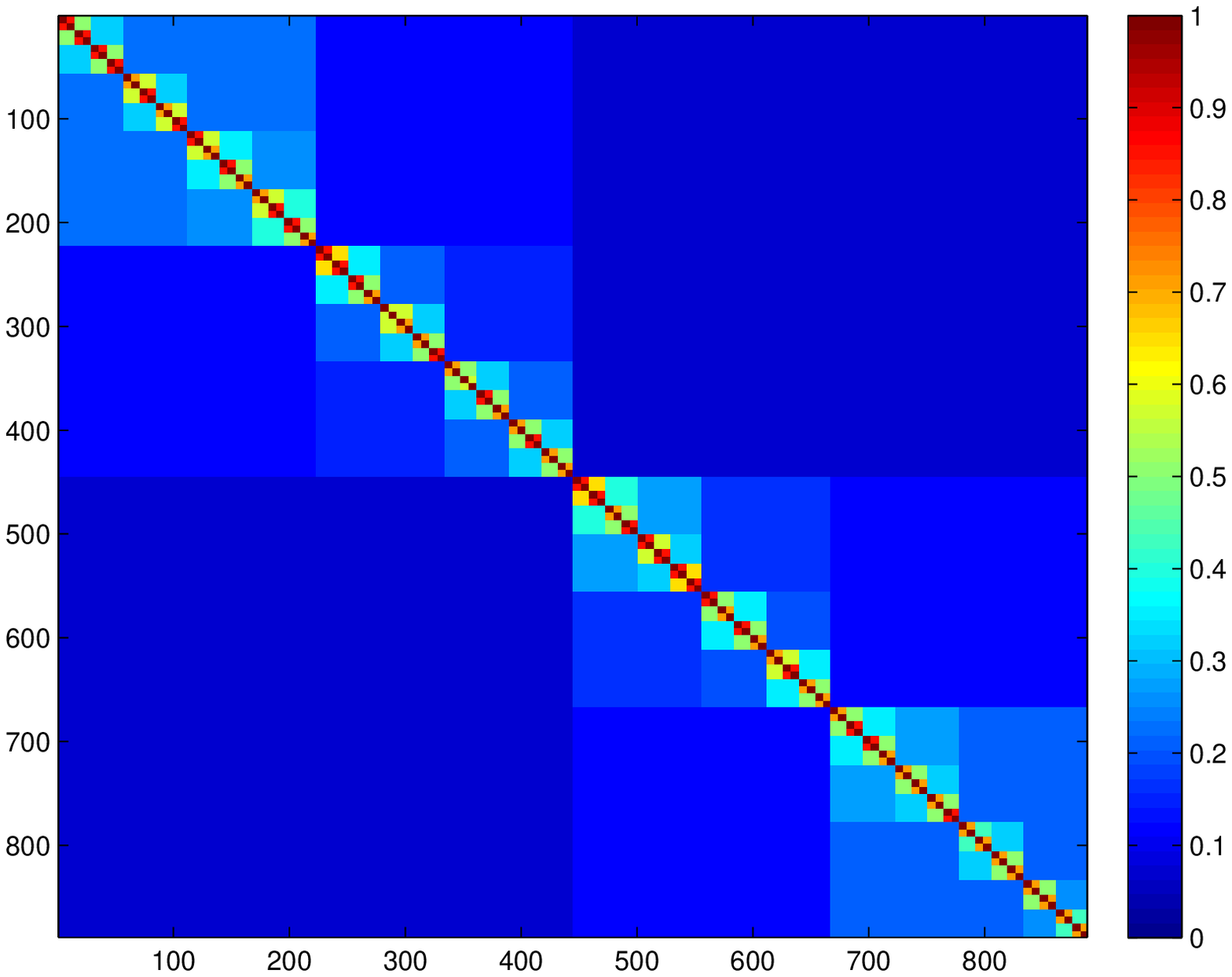} \\
\includegraphics[scale=0.35]{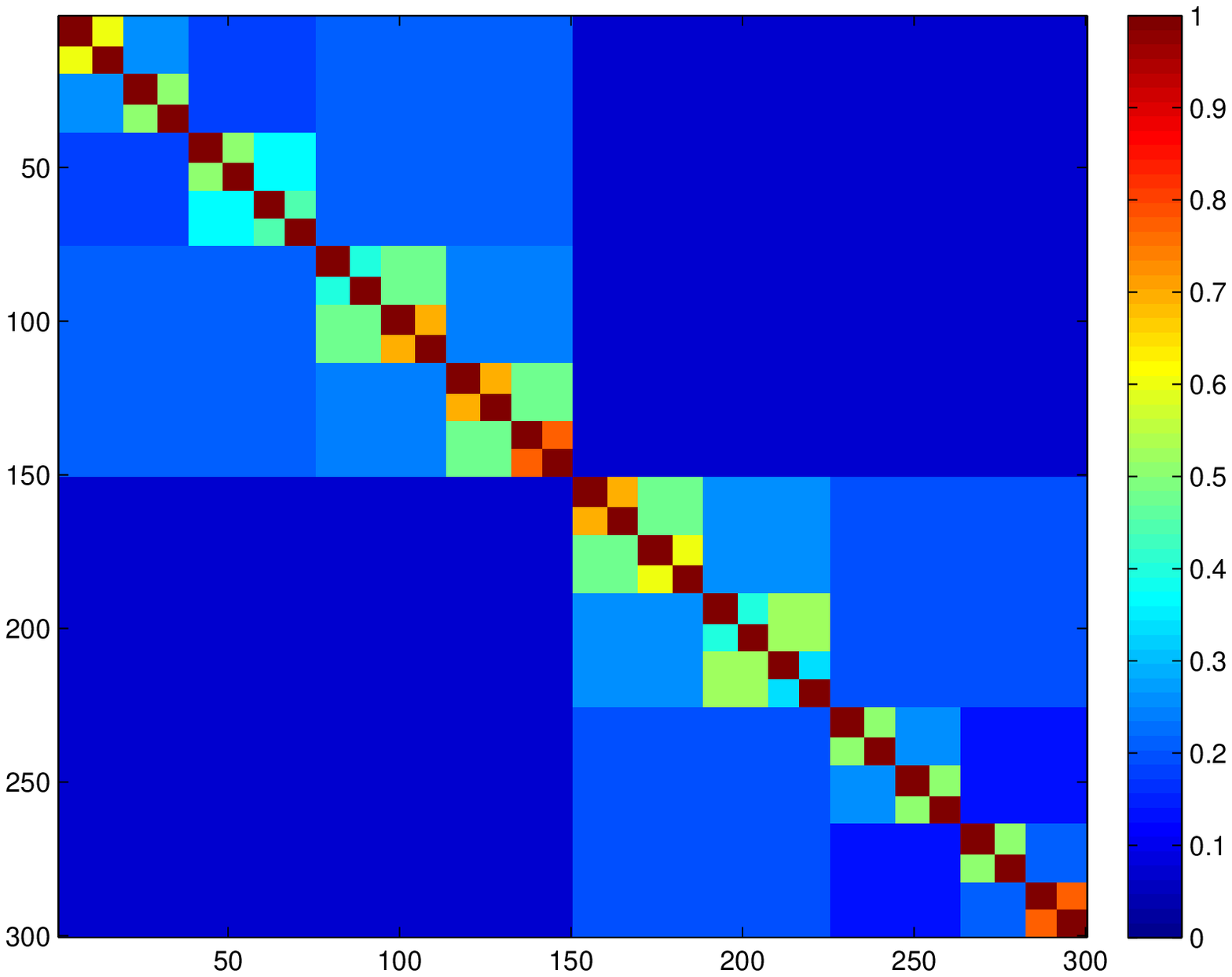}
\includegraphics[scale=0.35]{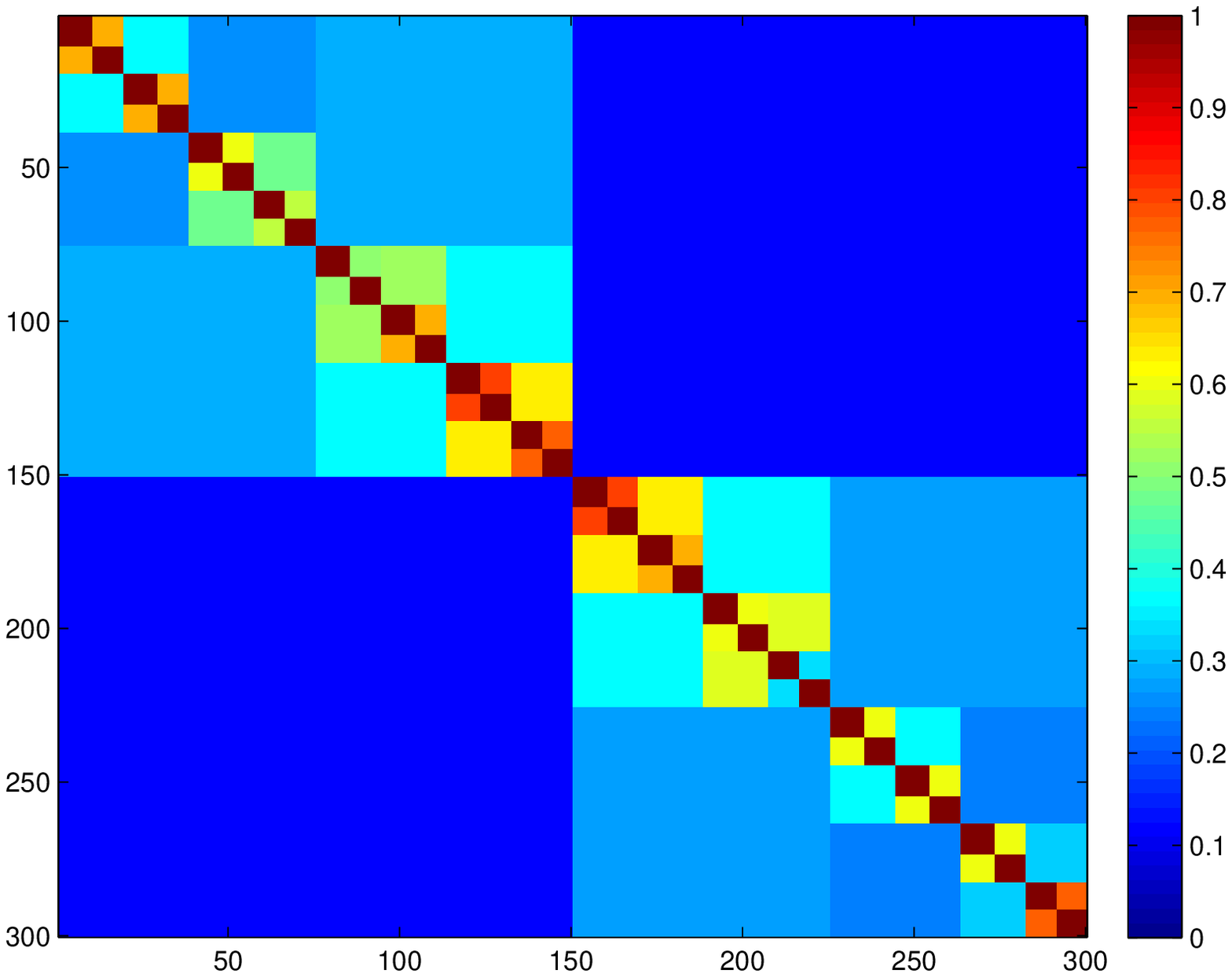} \\
\includegraphics[scale=0.35]{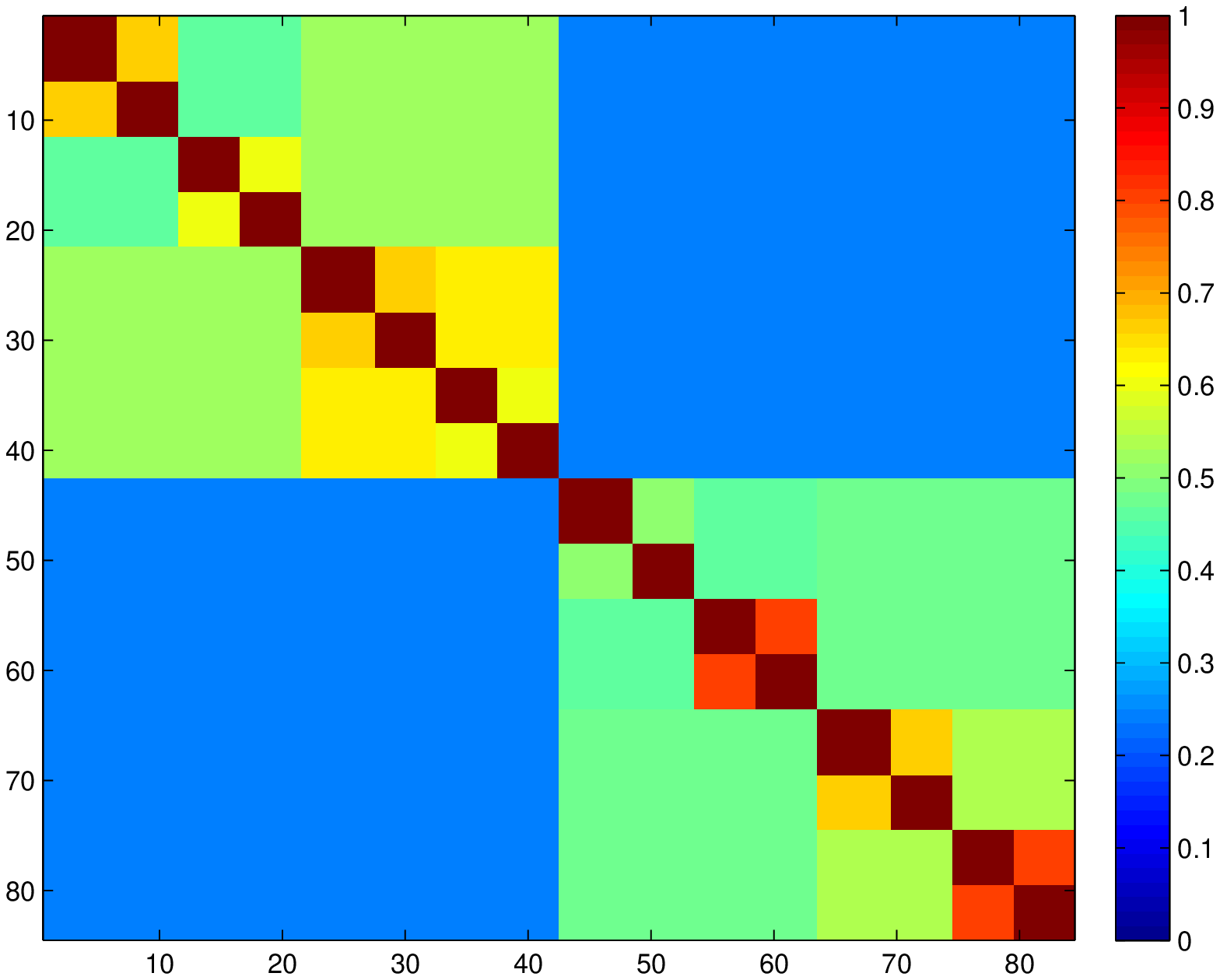}
\includegraphics[scale=0.35]{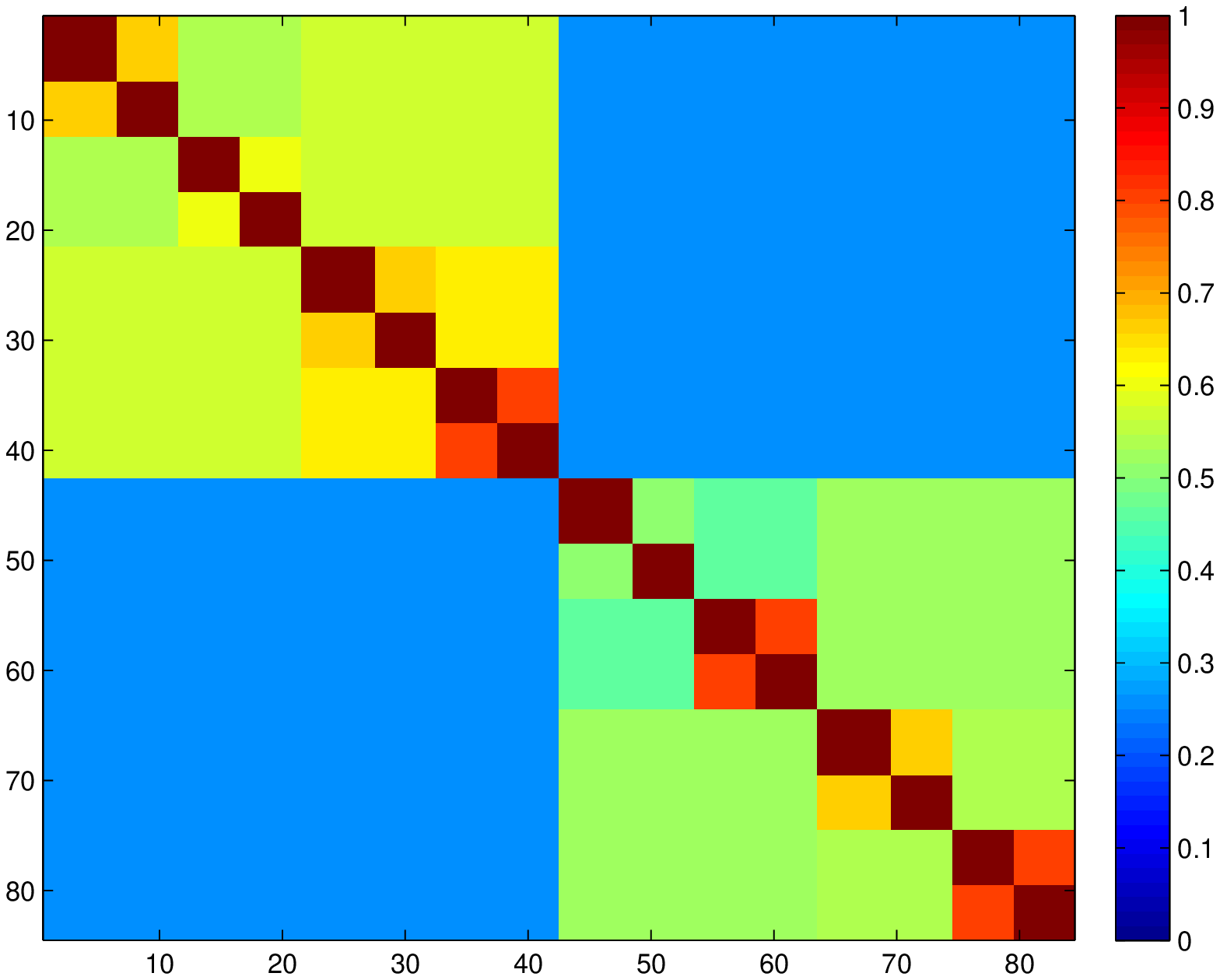}
\end{center}
\caption{Rank-structure of Schur complements of the Laplace operator (left column) as opposed to those of the Helmholtz operator (right column). From top to bottom, finer and finer tree-levels are considered. The colors indicate the relative rank of the corresponding sub-blocks. }\label{fig:5}
\end{figure}

In order to assess the validity of our theoretical framework, we analyze the rank-structure of the Schur complements that emerge from discretizations of the Laplace and the Helmholtz operators. A qualitative study that compares the Schur complements of the two operators is illustrated in Figure \ref{fig:5}. For a fixed problem size, we descend the tree focusing on the Schur complements relative to the first box on each level. At the top level, the Schur complement of the Laplace operator exhibits off-diagonal blocks that are rank-deficient, and more so than their counterparts relative to the Helmholtz operator. As we descend the tree, the difference between the Schur complements of the two operators becomes less evident. Further numerical experiments indicate that this behavior is consistent across problems of different sizes.

\begin{figure}[t]
\begin{center}
\includegraphics[scale=0.5]{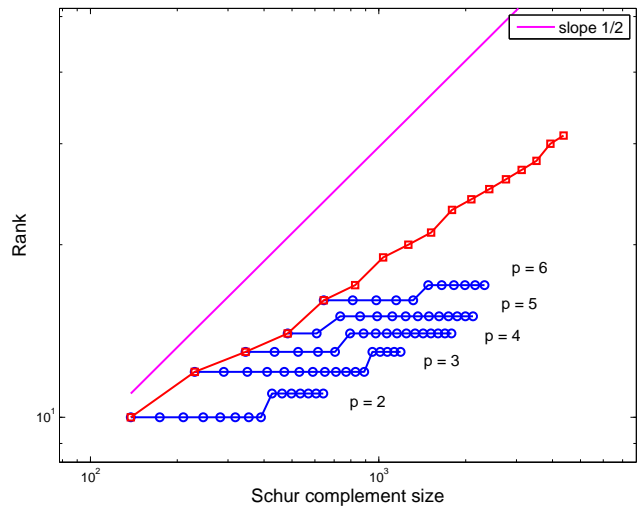}
\end{center}
\caption{Rank growth of Schur complements of Laplace operator. Data refer to ranks of largest off-diagonal block of complements on the top tree level. The problem size is increased through $h$-refinements (circles) or $p$-enrichments (squares.) \label{fig:9}}
\end{figure}

We proceed by studying the behavior of the HSS-rank $k$ of the Schur complements relative to the Laplace operator as the problem size $N$ is increased through $h$-refinements or $p$-enrichments, see Figure \ref{fig:9}. The numerical experiments showed a uniformity of behaviors of the Schur complements across all tree-levels. Thus, we limit our presentation to the Schur complement relative to the first box on the top level. When $h$-refinements are performed, $k$ is independent of the problem size while, in the case of $p$-enrichments, $k$ grows roughly as the square root of the size $n$ of the Schur complement. In the case of the Helmholtz operator, we postulate that $k$ grows as $\log n$ in the case of $h$-refinements, and as $n^{1/2}$ in the case of $p$-enrichments.

Taking into account that estimate (\ref{eq:42}) reduces to $\# B^\sigma = \mathrm{O}(N^{1/2})$ in the case of $h$-refinements, and to $\# B^\sigma = \mathrm{O}(N)$ in the case of $p$-enrichments, we can revisit the cost estimate (\ref{eq:38}) as follows:
\begin{description}
\item[{\it Laplace operator:}]
\begin{subequations}\label{eq:45}
\begin{align}
\text{cost of processing $\sigma$ on level $\ell$ ($h$-refinements)} & = \frac{1}{2^\ell}\,\mathrm{O}\big(N^{1/2}) \label{eq:39} \\
\text{cost of processing $\sigma$ on level $\ell$ ($p$-enrichments)} & = \frac{1}{2^\ell}\,\mathrm{O}\big(N^2) \label{eq:40}
\end{align}
\end{subequations}
\item[{\it Helmholtz operator:}]
\begin{subequations}\label{eq:55}
\begin{align}
\text{cost of processing $\sigma$ on level $\ell$ ($h$-refinements)} & = \frac{1}{2^\ell}\,\mathrm{O}\big(N^{1/2} \log^2 N) \label{eq:56} \\
\text{cost of processing $\sigma$ on level $\ell$ ($p$-enrichments)} & = \frac{1}{2^\ell}\,\mathrm{O}\big(N^2) \label{eq:57}
\end{align}
\end{subequations}
\end{description}
The estimates corresponding to the two scenarios, i.e., $h$-refinements as opposed to $p$-enrichments, are dramatically different. Apart from the different growth rates of $k$, this is dictated by the following geometrical argument. Given a fixed box, e.g., a box on the top tree-level, when $h$-refinements are performed, we observe a ``thinning'' of the boundary, namely the number of boundary dof's grows slower than the total number of dof's in the box. On the other hand, in the case of $p$-enrichments which corresponds to the limit $d \uparrow \infty$, the box contains a fixed number of elements, and no ``thinning'' of the boundary occurs. In applications of practical interest, the problem size is increased through adaptive strategies, that combine $h$-refinements and $p$-enrichments to produce optimal meshes. Estimates (\ref{eq:40}) and (\ref{eq:57}) should be regarded as the worst-case scenarios, which are not to be encountered in practice.

Although ``time'' (construction and application) is the ultimate measure of performance of a preconditioner, it is spectacularly implementation-dependent, and it makes little sense in the context of a non-optimized implementation. We take the following approach. Since the Schur complements are compressed by recursively performing ID's of appropriately select sub-blocks, we evaluate the total compression cost through estimate (\ref{eq:43}). More precisely, we select the anticipated rank $l$ according to the above discussion, and employ $l$ and the actual rank $k$, as returned by the ID, in the computation of the final estimate. We investigate the agreement between such empirical cost and the theoretical estimates (\ref{eq:45}) and (\ref{eq:55}). We measure the performance of the preconditioner in terms of GMRES iterations. If $l$ is properly selected, as a function of $N$, we expect the performance of the preconditioner not to degrade as the problem size is increased. We compare the performance of our preconditioner to that of a standard ILU preconditioner.

\begin{figure}[t]
\begin{center}
\subfigure[Problem size is increased through $h$-refinements.]{
\includegraphics[scale=0.39]{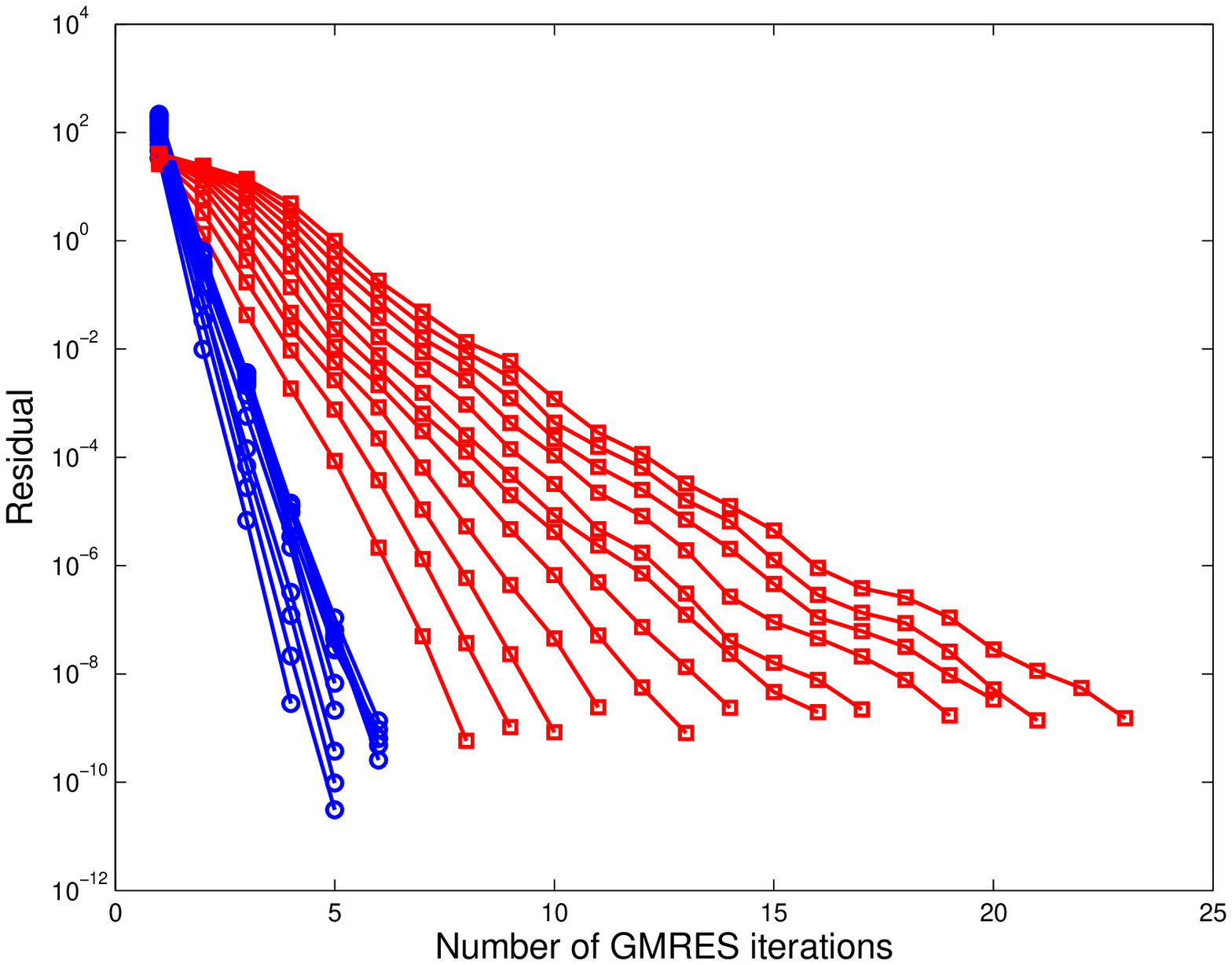}
\includegraphics[scale=0.39]{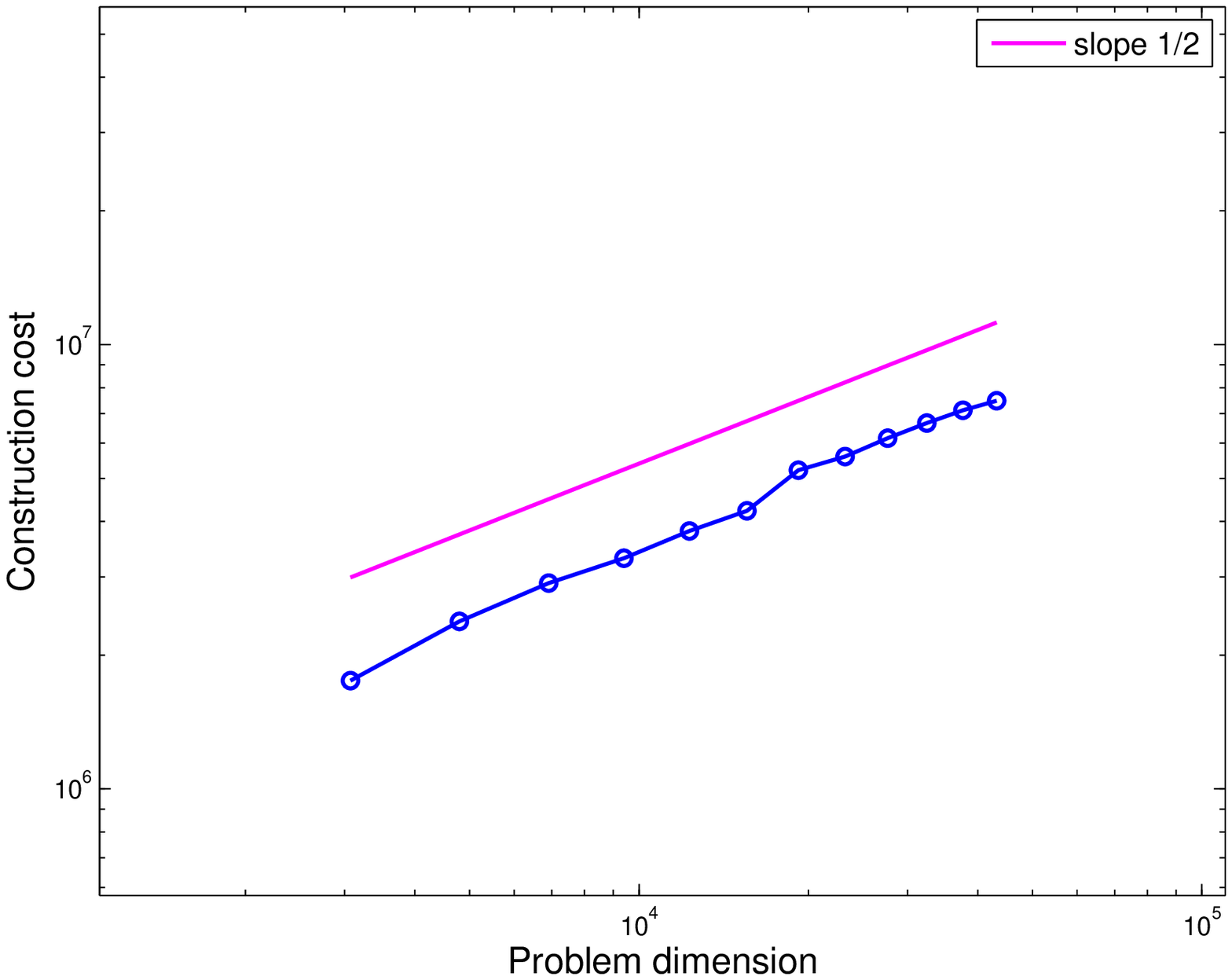}
\label{fig:72} } \\
\subfigure[Problem size is increased through $p$-enrichments.]{
\includegraphics[scale=0.39]{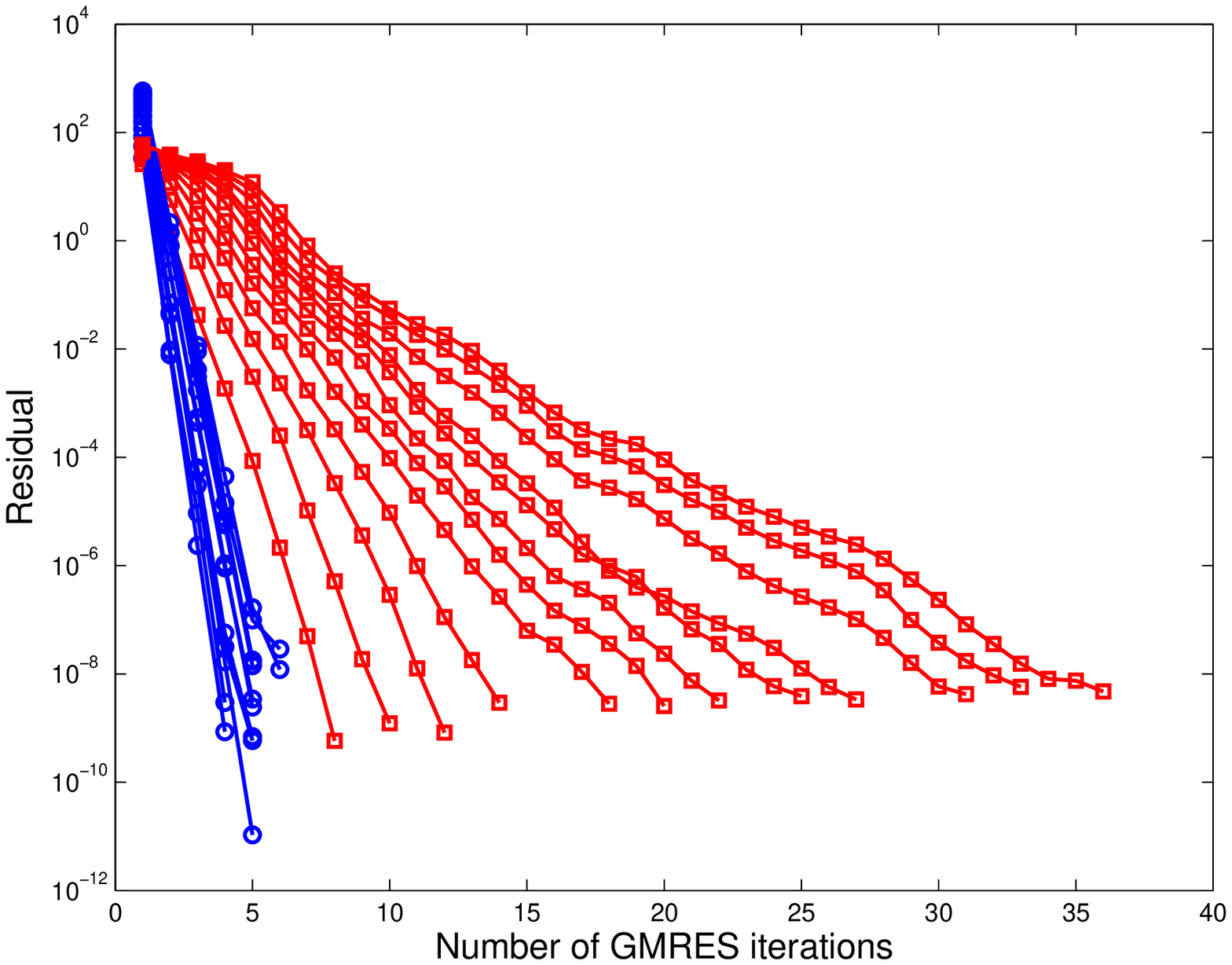}
\includegraphics[scale=0.39]{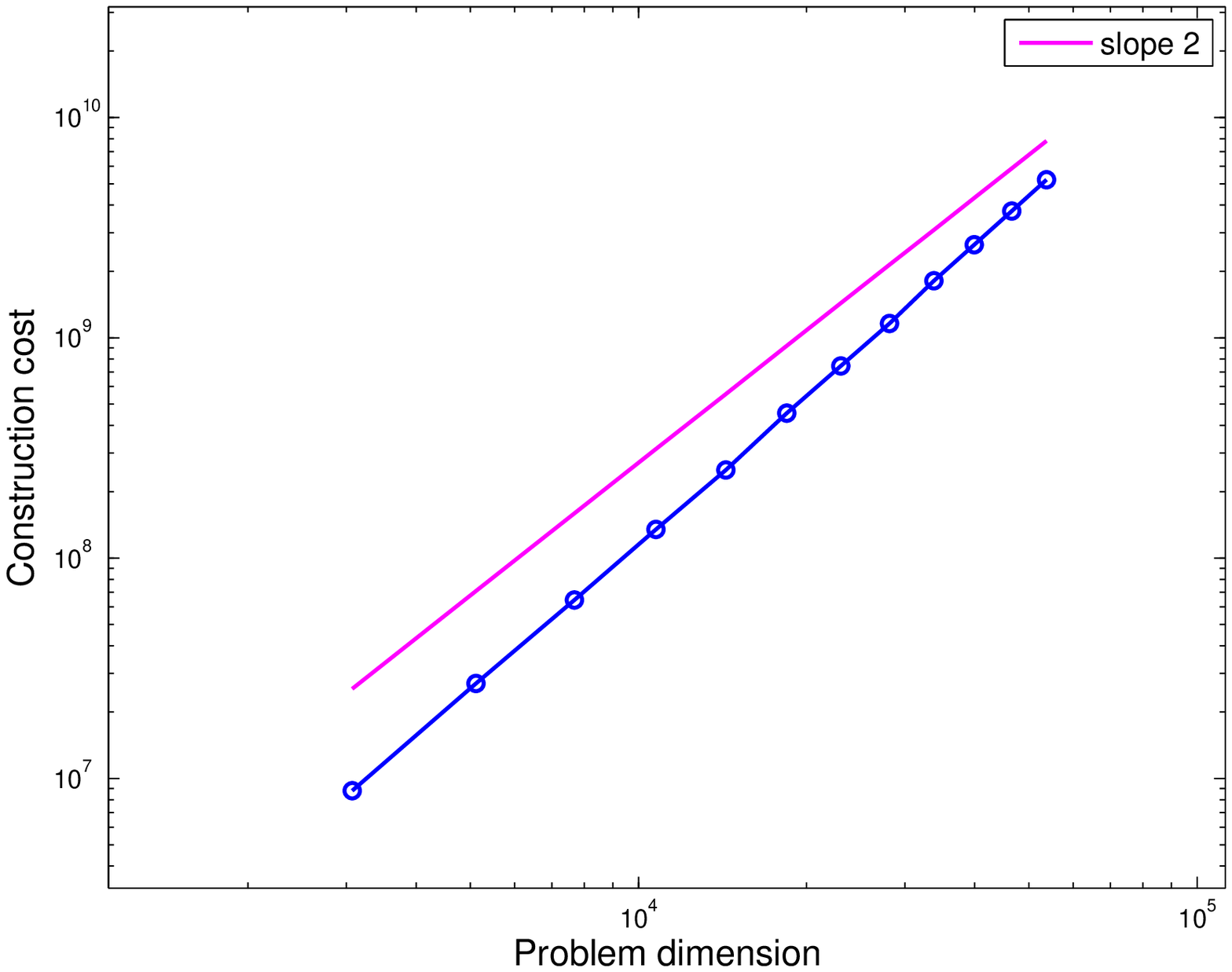}
\label{fig:74} }
\end{center}
\caption{Solution of a Poisson problem through GMRES. Circles refer to the proposed low-rank preconditioner, and squares refer to a standard ILU preconditioner. The performance of the low-rank preconditioner, i.e., the number of GMRES iterations, does not deteriorate as the problem size increases.}\label{fig:7}
\end{figure}

As a first example, we consider the Poisson equation on $\Omega$, with a homogenous Dirichlet boundary condition:
\begin{alignat*}{3}
-\Delta u &= f \qquad& &\text{in }\Omega \\
u & = 0	&& \text{on }\partial \Omega
\end{alignat*}
As anticipated, the problem is discretized using an Interior Penalty DG formulation. We build the preconditioner by recursively partitioning the domain into boxes, see Figure \ref{fig:5.1}, and study its behavior by comparing it to that of a standard ILU preconditioner. The construction cost agrees with the theoretical estimates (\ref{eq:45}), and the performance, i.e., the number of GMRES iterations, is independent of $h$, and $p$, see Figure \ref{fig:7}.

\begin{figure}[t]
\begin{center}
\subfigure[$c_0 = 10^{-4}$]{
\includegraphics[scale=0.38]{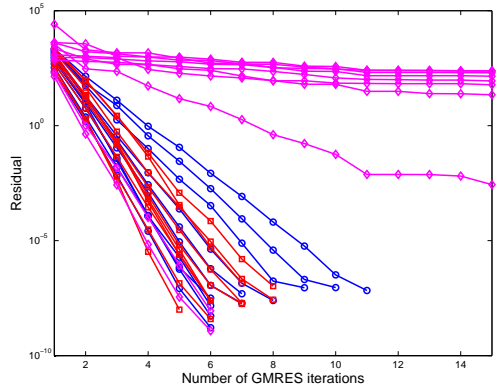}
}\subfigure[$c_0 = 1$]{
\includegraphics[scale=0.38]{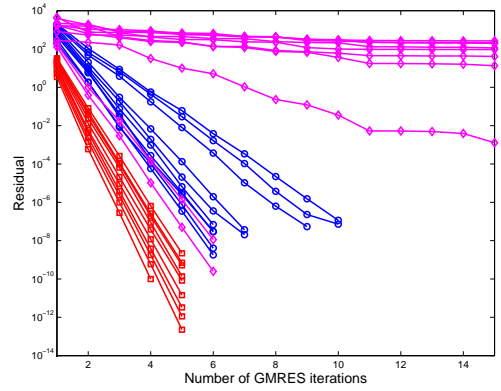}
}

\subfigure[$c_0 = 10^{2}$]{
\includegraphics[scale=0.38]{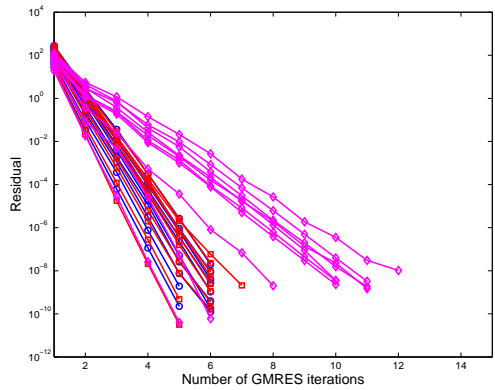}
}\subfigure[$c_0 = 10^{3}$]{
\includegraphics[scale=0.38]{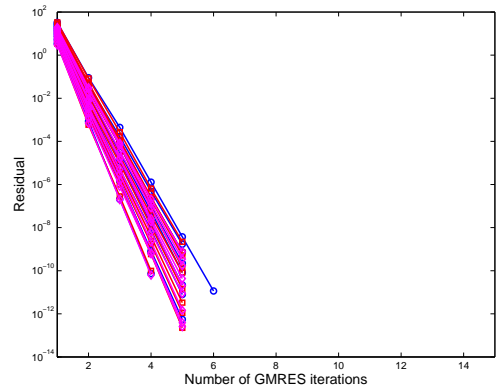}
}

\end{center}
\caption{Solution of a strongly anisotropic, i.e., $a_0 = 10^{-6}$, reaction/diffusion problem through GMRES. The problem size is increased through $h$-refinements. Circles refer to box partitioning of the domain, squares refer to horizontal slab partitioning, and diamonds refer to vertical slab partitioning. As $c_0$ increases and the problem transitions for diffusion-dominated to reaction-dominated, the performance of the low-rank preconditioner becomes independent of the partitioning scheme.}\label{fig:10}
\end{figure}

In order to investigate the role played by the partitioning scheme, we employ a constant coefficients reaction/diffusion problem, with an anisotropic diffusion tensor $A = \diag(a_{11} , a_{22})$:
\begin{subequations}
\begin{alignat}{3}
- \dive (A \, \nabla u) + c\,u &= f \qquad& &\text{in }\Omega  \label{eq:60} \\
u & = 0	&& \text{on }\partial \Omega
\end{alignat}
\end{subequations}
We construct the diffusion tensor so that it displays a strong anisotropy, i.e., $a_{11} \gg a_{22}$, and investigate the problem as it transitions from a diffusion-dominated regime to reaction-dominated regime. We define $a_0 = a_{22}/a_{11}$, $c_0 = c/a_{11}$, and $f_0 = f/a_{11}$, and proceed to normalize equation (\ref{eq:60}) by $a_{11}$:
\[
- \frac{\partial^2 u}{\partial x^2} - a_0\,\frac{\partial^2 u}{\partial y^2} + c_0\, u = f_0 \qquad \text{in } \Omega
\]
We set $a_0 \ll 1$, i.e., the diffusion is much stronger along the $x$-axis, as compared to the $y$-axis, and study the behavior of the preconditioner for $c_0 \ll 1$, $c_0 = 1$, and $c_0 \gg 1$, see Figure \ref{fig:10}. When the problem is not reaction-dominated, i.e., $c_0 \ll 1$ or $c_0 = 1$, the choice of the partitioning scheme is crucial. As in accordance with the nature of the problem, vertical slabs work very poorly, while horizontal slabs work excellently. The performance of the box partitioning scheme is comparable to that of the vertical slab partitioning, thus it can be regarded as a general purpose scheme. As the problem transitions to a reaction-dominated regime, the domain partitioning scheme becomes less and less relevant.

\begin{figure}[t]
\begin{center}
\subfigure[Problem size is increased through $h$-refinements.]{
\includegraphics[scale=0.39]{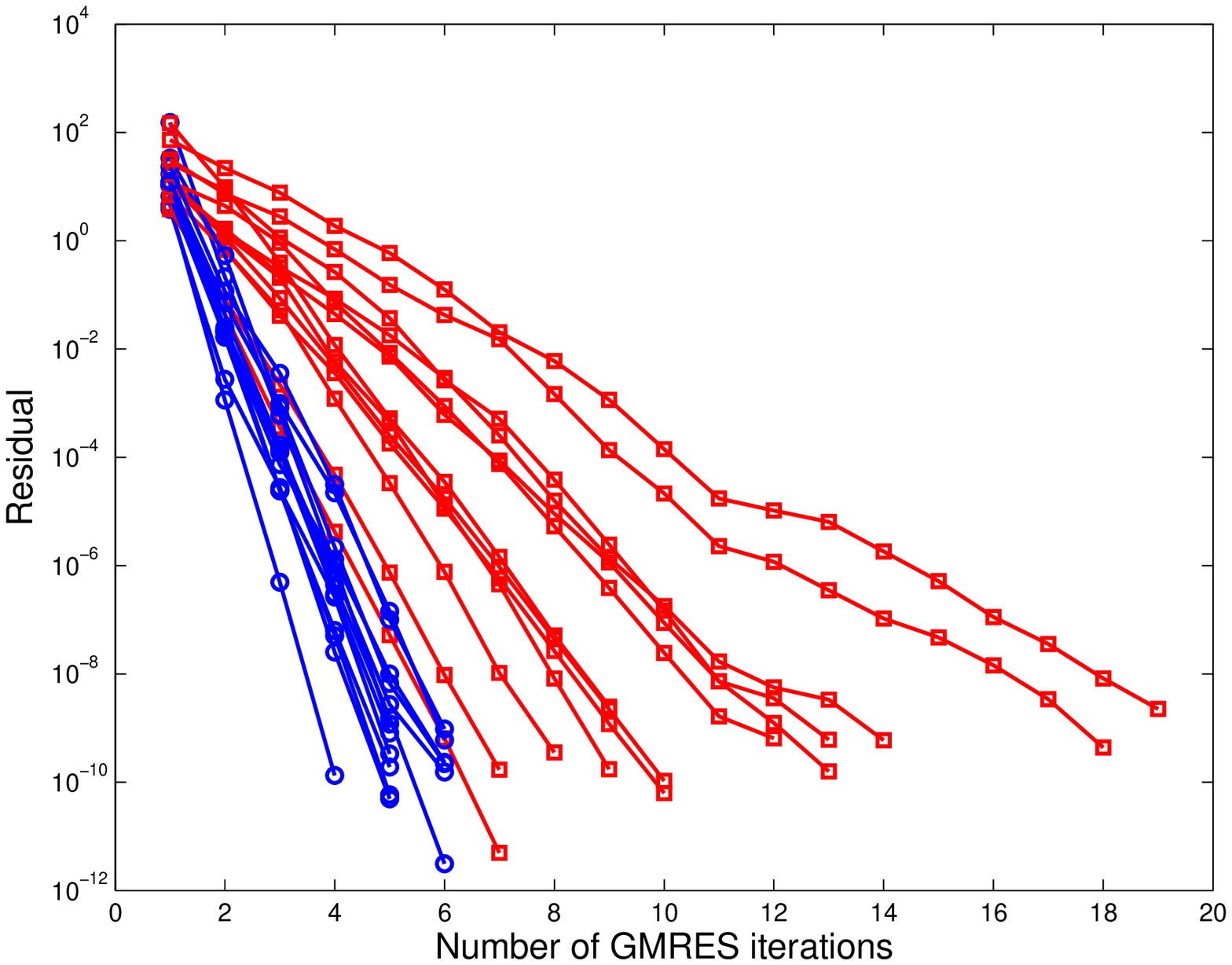}
\includegraphics[scale=0.39]{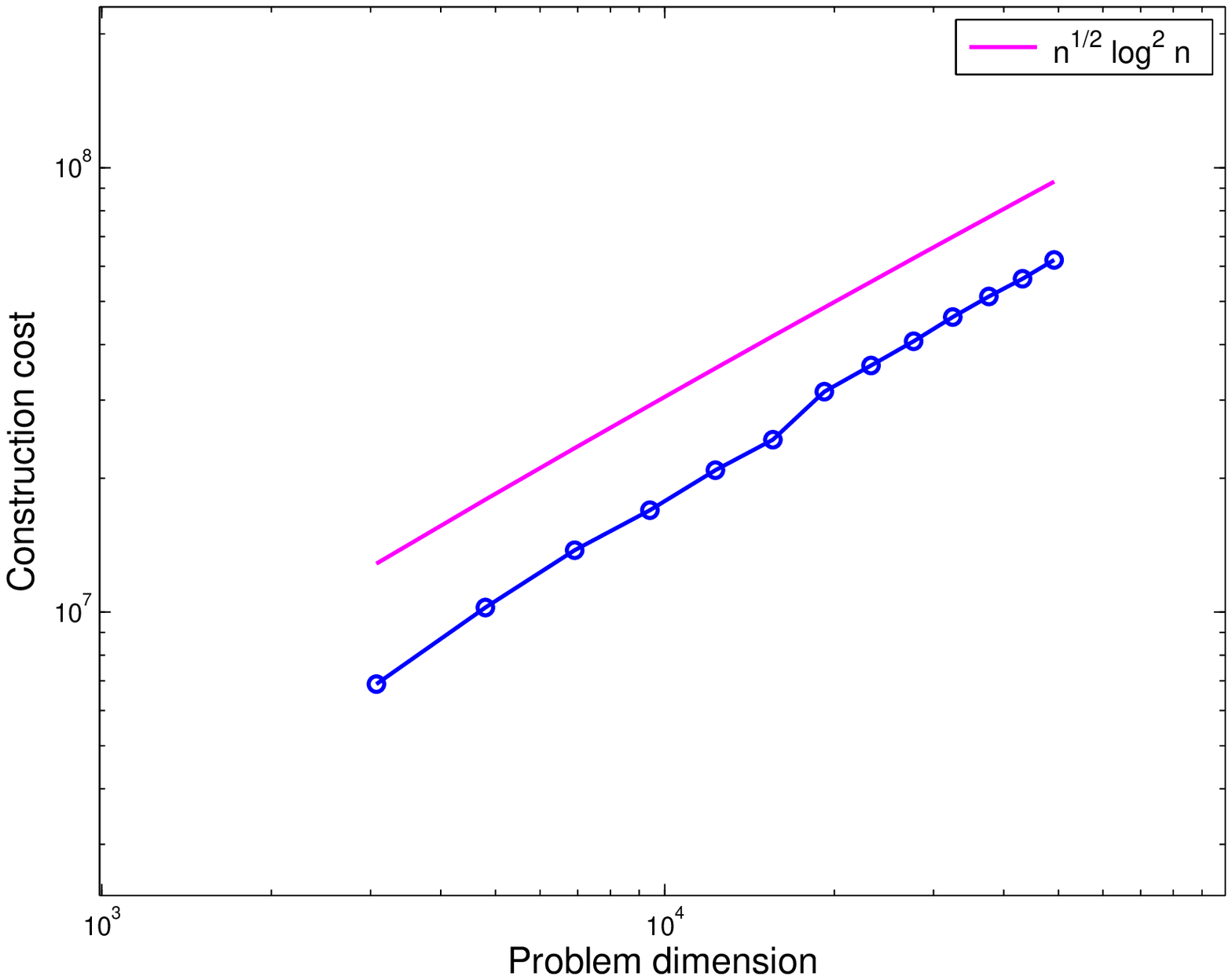}
} \\
\subfigure[Problem size is increased through $p$-enrichments.]{
\includegraphics[scale=0.39]{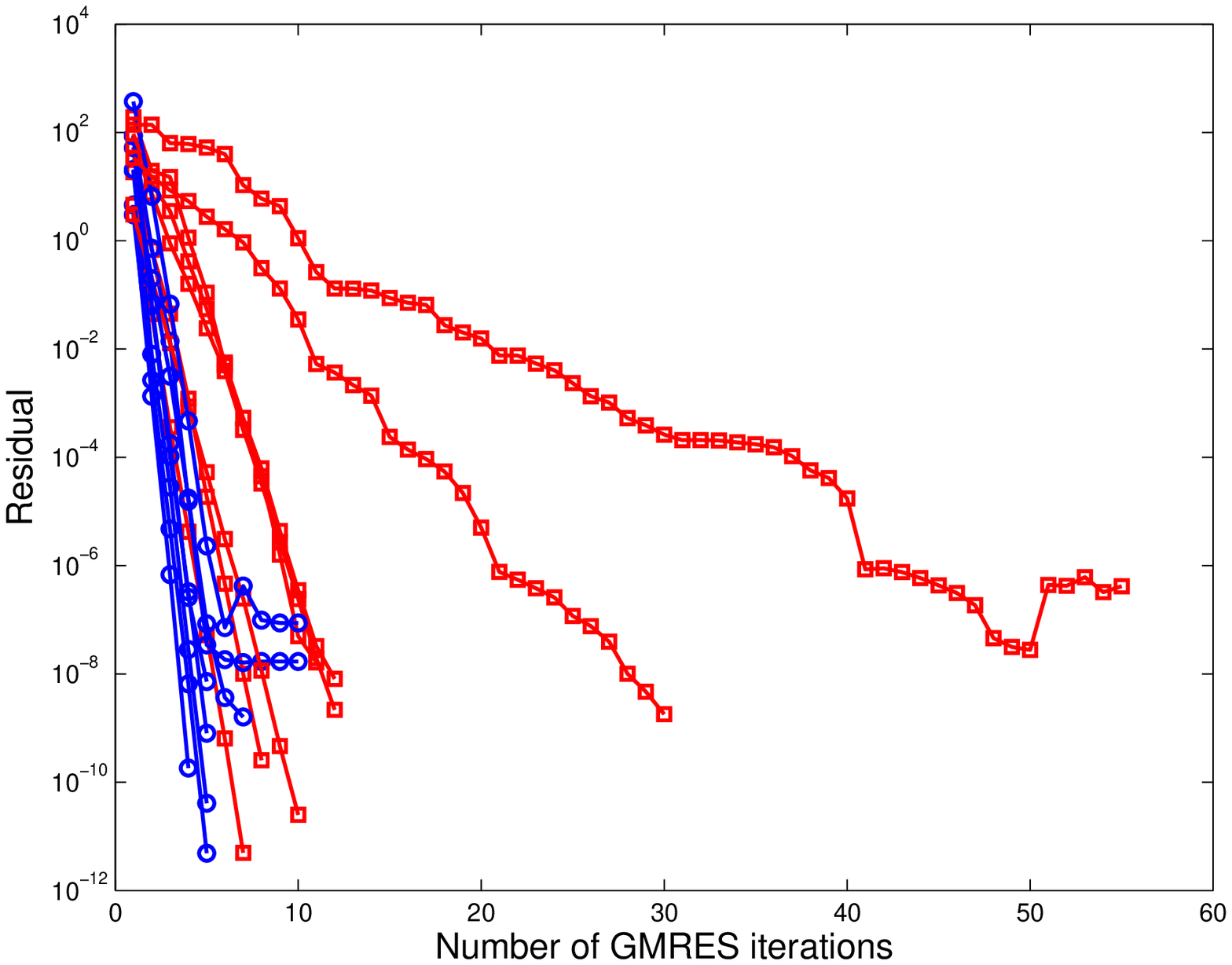}
\includegraphics[scale=0.38]{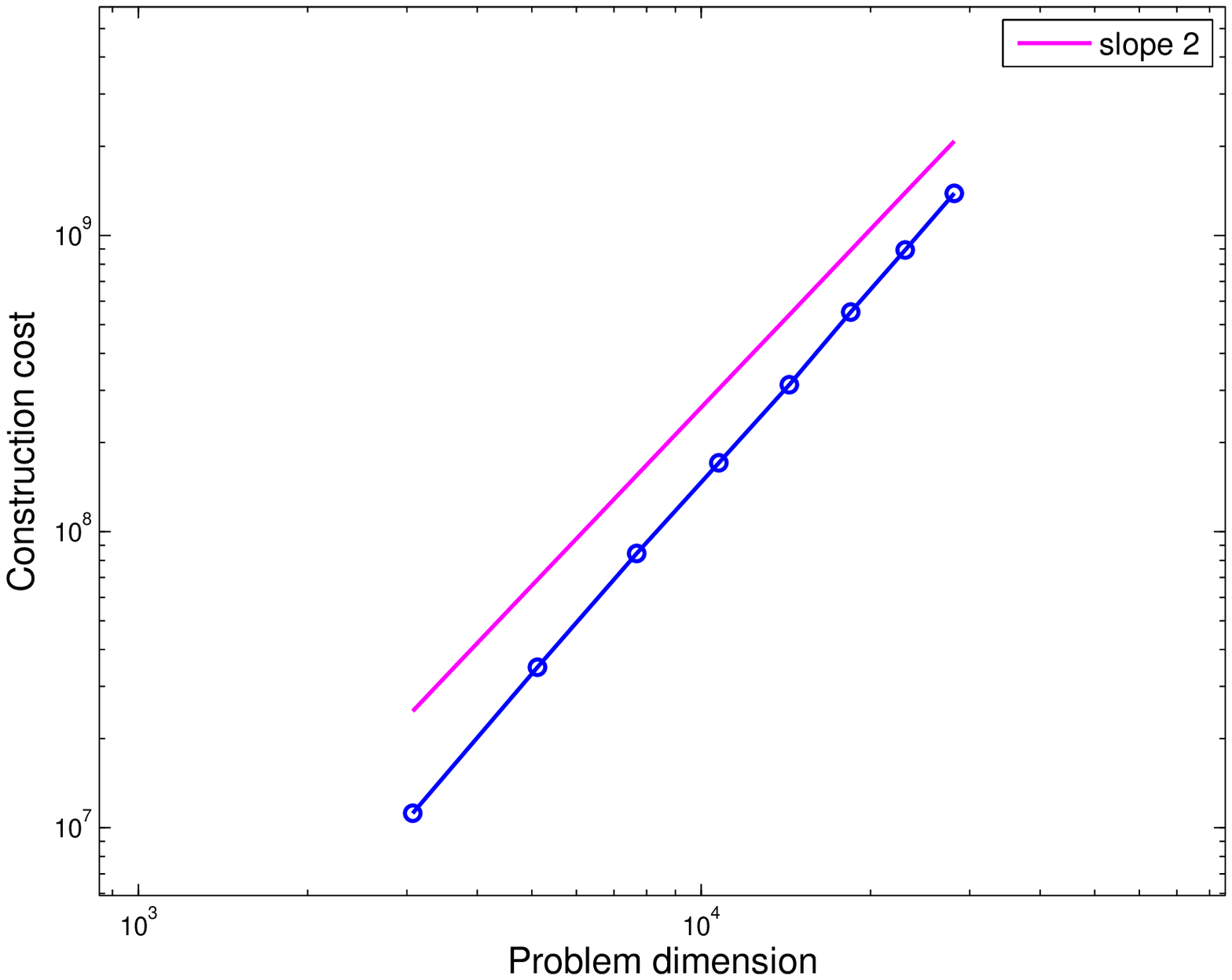}
}
\end{center}
\caption{Solution of a Helmholtz problem through GMRES. Circles refer to the proposed low-rank preconditioner, and squares refer to a standard ILU preconditioner. The performance of the low-rank preconditioner, i.e., the number of GMRES iterations, does not deteriorate as the problem size increases.}\label{fig:8}
\end{figure}

Finally, we move to indefinite problems. Let us consider the Helmholtz equation on $\Omega$ with a homogenous Dirichlet boundary condition:
\begin{alignat*}{3}
-\Delta u - \kappa^2 u &= f \qquad& &\text{in }\Omega \\
u & = 0	&& \text{on }\partial \Omega
\end{alignat*}
This boundary value problem describes propagation of forced waves inside a soft cavity. Although exterior scattering problems are often of greater interest, the cavity problem is computationally more challenging because of possible resonances.  Remarkably, the behavior of the preconditioner is qualitatively similar to that observed in the case of the Poisson problem. The construction cost is in agreement with estimates (\ref{eq:55}) and the performance of the preconditioner does not deteriorate as the problem size is increased through $h$-refinements or $p$-enrichments, see Figure \ref{fig:8}. 

\begin{figure}[t]

\begin{center}

\subfigure[Domain partitioning schemes (enclosure $D$ is indicated in gray.)]{
\begin{tikzpicture}[scale = 0.05]

\draw[black] (0,0) rectangle (80,80);

\path[fill=gray] (10,30) rectangle (70,40);
\draw[line width = 2.5](1,40) -- (79,40);

\draw[line width = 1.5](40,1) -- (40,38);
\draw[line width = 1.5](40,42) -- (40,79);

\end{tikzpicture}
\hspace{1cm}
\begin{tikzpicture}[scale = 0.05]
\path[fill=gray] (10,30) rectangle (70,40);

\draw[black] (0,0) rectangle (80,80);
\draw[line width = 2.5](1,40) -- (79,40);

\draw[line width = 1.5](1,20) -- (79,20);
\draw[line width = 1.5](1,60) -- (79,60);

\end{tikzpicture}
\hspace{1cm}
\begin{tikzpicture}[scale = 0.05]

\path[fill=gray] (10,30) rectangle (70,40);

\draw[black] (0,0) rectangle (80,80);

\draw[line width = 2.5](40,1) -- (40,79);

\draw[line width = 1.5](20,1) -- (20,79);
\draw[line width = 1.5](60,1) -- (60,79);

\end{tikzpicture}\label{fig:9.1}
}

\subfigure[Circles refer to box partitioning, squares refer to horizontal slab partitioning, diamonds refer to vertical slab partitioning]{
\includegraphics[scale=0.5]{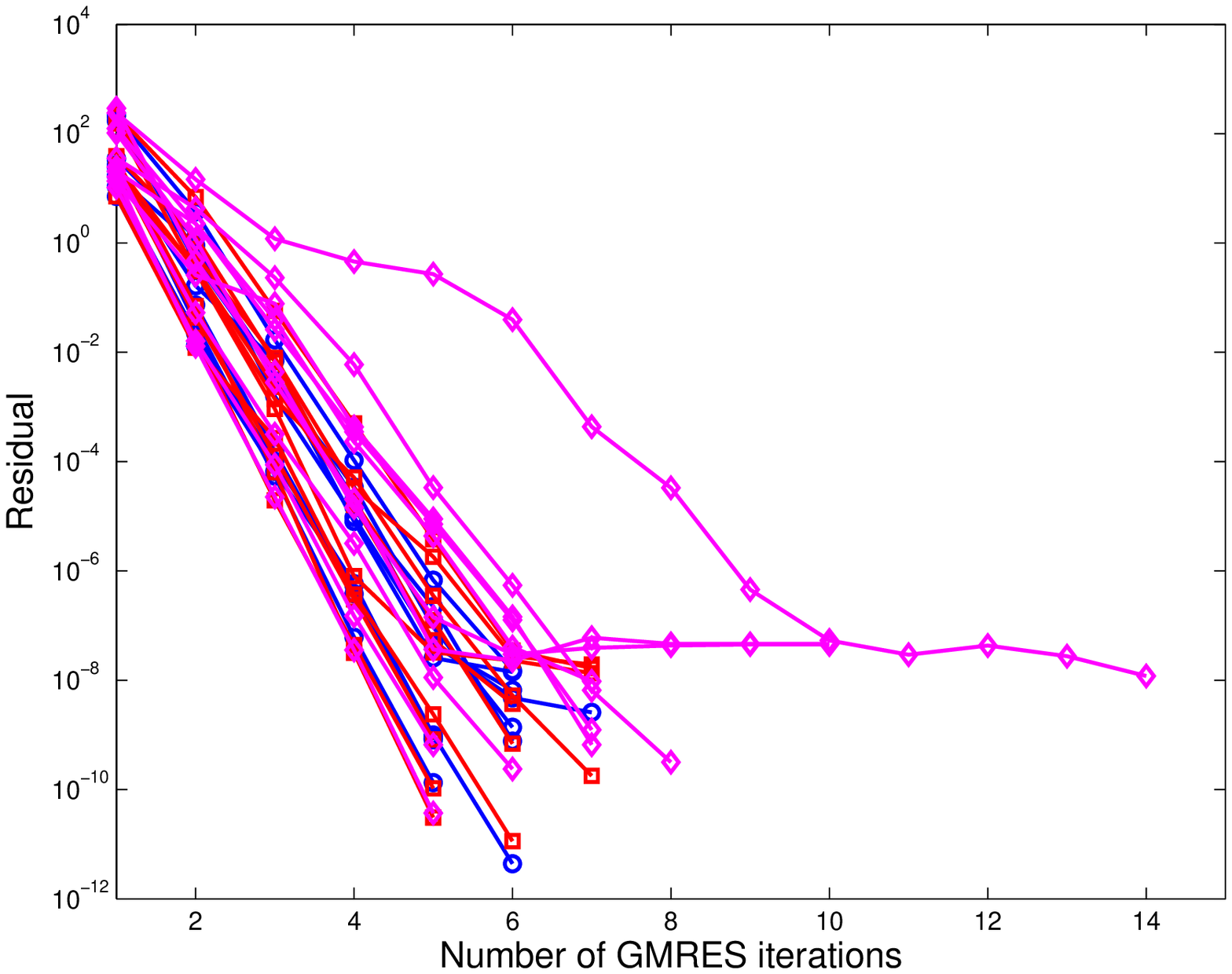} \label{fig:9.2}
}
\end{center}
\caption{Solution of a high-material-contrast Helmholtz problem through GMRES. The material enclosure $D$ has non-dimensional density $\varrho_D = 10^3$, while the remaining portion of the domain, i.e., $\Omega \setminus D$, has density $\varrho_0 = 1$. The problem size is increased through $h$-refinements. The performance of the low-rank preconditioner exhibits a mild dependency on the partitioning scheme.}
\end{figure}

As a last example, we consider the following Helmholtz problem:
\begin{alignat*}{3}
-\dive\bigg( \frac{1}{\varrho} \, \nabla u\bigg) - \kappa^2 u &= f \qquad& &\text{in }\Omega \\
u & = 0	&& \text{on }\partial \Omega
\end{alignat*}
This problem is of interest in, e.g., topology optimization, see \cite{bendsoe2009topology}. In a nutshell, topology optimization is an iterative method to create highly optimized designs by determining a distribution of material, namely $\varrho$, that fulfills a specific task, in a locally optimal manner. Consequently, the density $\varrho$ varies by orders of magnitude over the domain. For ease of implementation, we assume it to be constant over each element. We place a rectangular enclosure $D$ with density $\varrho_D$ inside the domain $\Omega$, and assume $\varrho_D \gg \varrho_0$, where $\varrho_0$ is the density in $\Omega \setminus D$, see Figure \ref{fig:9.1}. We observe a very mild dependency of the behavior of the preconditioner upon the domain partitioning scheme, see Figure \ref{fig:9.2}. As expected, the horizontal slab partitioning delivers the best performance, which is virtually matched by that of the box partitioning. Although the vertical slab partitioning \emph{a priori} appears as the least appropriate one, its performance is not dissimilar from those of the other partitioning schemes.

\section{Conclusions}\label{sec:5}
We have presented the construction of a preconditioner that exploits low-rank compression of Schur complements. The construction can be viewed as a variant of the well-known nested dissection algorithm, since it employs a reordering of the degrees of freedom which allows for an advantageous elimination order. Our approach, originally inspired by the work of \cite{gillman2014n}, follows a black-box approach that gives the construction the flexibility to be applied to a number of discretization techniques, such as finite differences, CG finite elements and DG finite elements. The preconditioner can be applied within linear complexity, and we provide an estimate of the construction cost. Such estimate depends widely upon the strategy used to vary the problem size ($h$-refinement as opposed to $p$-enrichments, in the case of finite elements approximations), with a worst-case-scenario of quadratic growth. Although this issue requires further investigation, we believe that for applications of practical interest, the construction cost is within linear growth. We tested the performance of the preconditioner on DG approximations of elliptic as well as hyperbolic problems, see \cite{hesthaven2007nodal}, and demonstrated its robustness. More specifically, the preconditioned system is solved within a number of GMRES iterations that is independent of both the mesh size and the order of approximation. The choice of DG finite elements approximations was dictated by implementation convenience, namely the reordering of the degrees of freedom has a straightforward geometrical interpretation, and should not be viewed as a limitation of the applicability of the preconditioner. In principle, the reordering of the degrees of freedom can be performed with any graph-partitioning software. 

The construction of preconditioners for linear systems that arise from wave propagation phenomena is notoriously challenging, see \cite{toselli2005domain}. The proposed preconditioner is based on a completely general construction and has proven effective for a number of problems. In terms of future research developments, the most pressing issue is to verify the agreement between the actual construction cost and its theoretical estimate through an optimized implementation. This will allow us to assess whether the asymptotic region is reached for problem sizes of practical relevance. Finally, we should investigate the effectiveness of the preconditioner over a larger class of problems, e.g., coupled multi-physics problems, and explore the possibility of adapting the construction to time-dependent problems.

\section*{References}
\bibliography{paper}
\bibliographystyle{plain}

\end{document}